\documentclass[12pt]{amsart}
\usepackage[dvips]{graphicx}
\usepackage{color}
\input epsf
\def\cal{\mathcal}

\def\RR{\mathbb R}
\def\HH{\mathbb H}
\def\ts{{\cal T}(S)}
\def\a{\alpha}
\def\s{\sigma}
\def\L{{\cal L}}
\def\M{{\cal M}}
\def\T{{\cal T}}

\def\tthick{{\cal T}_{thick}(S)}
\def\da{d_\a(\nu^+,\nu^-)}
\def\dd{\partial}
\def\ep{\epsilon} 

\def\np{{\nu^+}}
\def\nm{{\nu^-}}
\def\ua{{\underline \a}}

\def\Teich{Teichm\"uller }

\def\mul{\stackrel{{}_\ast}{\asymp}}

\def\b{\beta}

\def\g{\gamma}
\def\t{\tau}

\def\E{{\cal E}}

\def\G{{\cal G}}
\def\M{{\cal M}}
\def\Q{{\cal Q}}
\def\S{{\cal S}}

\def\teich{{\cal T}}
 
 \def\ML{\mathop{\cal{ ML}} (S)}
 \def\mod{\mathop{\rm{Mod}}}
   
   \def\Ext{\mathop{{Ext}}}

\newtheorem{theorem}{Theorem}[section]
\newtheorem{definition}[theorem]{Definition}
\newtheorem{lemma}[theorem]{Lemma}
\newtheorem{proposition}[theorem]{Proposition}
\newtheorem{corollary}[theorem]{Corollary}

\newtheorem{introthm}{Theorem}

\begin{document}
\title{Lines of minima and Teichm\"uller geodesics}
\begin{abstract}
  For two measured laminations $\nu^+$ and $\nu^-$ that fill up a
  hyperbolizable surface $S$ and for $t \in (-\infty, \infty)$, let
  $\L_t$ be the unique hyperbolic surface that minimizes the length
  function $e^tl({\nu^+}) + e^{-t} l({\nu^-})$ on \Teich space.  We
  characterize the curves that are short in $\L_t$ and estimate their
  lengths. We find that the short curves coincide with the curves that
  are short in the surface $\G_t$ on the \Teich geodesic whose
  horizontal and vertical foliations are respectively, $e^t \np$ and
  $e^{-t} \nm$.  By deriving additional information about the twists
  of $\np$ and $\nm$ around the short curves, we estimate the \Teich
  distance between $\L_t$ and $\G_t$.  We deduce that this distance
  can be arbitrarily large, but that if $S$ is a once-punctured torus
  or four-times-punctured sphere, the distance is bounded
  independently of $t$.
\end{abstract}

\author{Young-Eun Choi} 
\email{choiye@psu.edu}

\author{Kasra Rafi}
\email{rafi@math.uconn.edu}

\author{Caroline Series}
\email{cms@maths.warwick.ac.uk}
\maketitle
\date{}
\section{Introduction}
\label{sec:introduction}

Suppose that $\nu^+$ and $\nu^-$ are measured laminations which fill
up a hyperbolizable surface $S$.  The object of this paper is to
compare two paths in the Teichm\"uller space $\T(S)$ of $S$ determined
by $\nu^+$ and $\nu^-$. The first is the Teichm\"uller geodesic $\G =
\G(\nu^+,\nu^-)$, whose time $t$ Riemann surface $\G_t \in \G$
supports a quadratic differential $q_t$ whose horizontal and vertical
foliations are $\nu^+_t = e^t\nu^+$ and $\nu^-_t =e^{-t}\nu^-$,
respectively \cite{G-M}.  The second is the Kerckhoff line of minima
$\L= \L(\nu^+,\nu^-)$ \cite{lminima}. At time $t$, $\L_t \in \cal L$
is the unique hyperbolic surface that minimizes the length function
$l(\nu^+_t) + l(\nu^-_t) = e^{t}l(\nu^+) + e^{-t}l(\nu^-)$ on $\ts$.
(Recall that $\G$ can be characterized as the locus in $\T(S)$ where
the product of the extremal lengths $\Ext(\np)\Ext(\nm)$ is minimized
\cite{G-M}.)  Lines of minima have many properties in common with
Teichm\"uller geodesics, see~\cite{lminima}, and have been shown to be
closely linked to deforming Fuchsian into quasi-Fuchsian groups by
bending, see~\cite{serieslm}.

We are interested in comparing the two trajectories $\G$ and $\L$, in
particular, to see whether or not they remain a bounded distance
apart. If both $\G_t$ and $\L_t$ are contained in the thick part of
$\ts$, it is relatively easy to show that the \Teich distance between
them is uniformly bounded independently of $t$, see
Theorem~\ref{thm:warmup}.  A more surprising fact is that in general,
the sets of short curves on $\G$ and $\L$ are the same.  Writing
$l_{\G_t}(\a), l_{\L_t}(\a) $ for the geodesic lengths of a simple
closed curve $\a$ in the hyperbolic metrics on $\G_t,\L_t$
respectively, we prove:
\begin{introthm}[Proposition~\ref{cor:shorteronL} and
  Corollary~\ref{cor:shortsame}]
\label{shortonboth}
The set of short curves on $\G_t$ and $\L_t$ coincide. More precisely,
there exist universal constants $\epsilon_1,\ldots, \epsilon_4>0$ such
that for each $t$, $ l_{\G_t}(\a) < \epsilon_1$ implies $l_{\L_t}(\a)
< \epsilon_2$ and $l_{\L_t}(\a) < \epsilon_3$ implies $l_{\G_t}(\a) <
\epsilon_4$.  \end{introthm}

Finding combinatorial estimates for the lengths of these short curves
occupies the main part of the paper and leads to a coarse estimate of
the distance between $\G_t$ and $\L_t$.  It turns out that along both
$\G$ and $\L$ there are two distinct reasons why a curve $\alpha$ can
become short: either the relative twisting of $\np$ and $\nm$ about
$\alpha$ is large, or $\nu^+$ and $\nu^-$ have large relative
complexity in $S \setminus \a$ (the completion of the surface $S$
minus $\alpha$), in the sense that every essential arc or closed curve
in $S \setminus \a$ must have large intersection with $\np$ or $\nm$.
Our results give sufficient control to construct examples which show
that $\G_t$ and $\L_t$ may or may not remain a bounded distance apart.

\medskip
 
The estimates for curves which become short along $\G$ are based on
Rafi~\cite{rafi1}.  For convenience we say a curve is `extremely
short' on a given surface if its hyperbolic length is less than some
fixed constant $\ep_0 > 0$ defined in terms of the Margulis constant,
see Section~\ref{sec:ttdecomp}.  Rafi's results imply:
\begin{introthm}[Theorem~\ref{prop:shortingeo}]
\label{shortong}
Suppose that $\a$ is extremely short on $\G_t$. Then
$$
\frac{1}{l_{\mathcal G_t}(\a)} \asymp \max \{ D_t(\a), \log K_t(\a) \}.$$
\end{introthm}
\noindent  The terms $D_t(\a)$ and $K_t(\a)$ correspond 
respectively to the relative twisting
and large relative complexity mentioned above.  More precisely,
$$D_t(\a) = e^{-2|t-t_{\a}|} d_{\a}(\np,\nm)$$
where $t_{\a}$ is the
\emph{balance time} at which $i(\a,\nu^+_t) = i(\a,\nu^-_t)$ and $
d_{\a}(\np,\nm)$ is the relative twisting, that is, the difference
between the twisting of $\np$ and $\nm$ around $\a$, see
Section~\ref{sec:reltwist}.  The term $K_t(\a)$ depends on the
(possibly coincident) thick components $Y_1,Y_2$ that are adjacent to $\a$
in the thick-thin decomposition of 
$\G_t$.  Let $q_t$ be the area $1$ quadratic
differential on $\G_t$ whose horizontal and vertical foliations are
respectively $\nu_t^+$ and $\nu_t^-$. Associated to $q_t$ is a
singular Euclidean metric; we denote the geodesic length of a curve
$\gamma$ in this metric by $l_{q_t}(\gamma)$, see
Section~\ref{sec:quaddiffs}.  By definition $$K_t(\a) = \max \Big\{
\frac{\lambda_{Y_1}}{l_{q_t}(\a)}, \frac{\lambda_{Y_2}}{l_{q_t}(\a)}
\Big\}
$$
where $\lambda_{Y_i}$ is the length of the shortest non-trivial,
non-peripheral simple closed curve on $Y_i$ with respect to the
$q_t$-metric, see Section~\ref{sec:thickthin}.
 
\medskip 

One of the main results of this paper is a similar characterization of
curves which become short along $\L(\np,\nm)$.  We prove that the
hyperbolic length $l_{\L_t}(\a)$ of a short curve in $\L_t$ is
estimated as follows:
\begin{introthm}[Theorem~\ref{cor:final}]
\label{shortonL}
Suppose that $\a$ is extremely short on $\L_t$. Then
$$\frac{1}{l_{\L_t}(\a)} \asymp \max \{ D_t(\a), \sqrt{K_t(\a)} \}.$$
\end{introthm}
The main tool in the proof is the well-known derivative formula of
Kerckhoff~\cite{kercktwist} and Wolpert~\cite{wolpert} for the
variation of length with respect to Fenchel-Nielsen twist, together
with the extension proved by Series~\cite{cmswolpert} for variation
with respect to the lengths of pants curves.

\medskip 

To estimate the Teichm\"uller distance between two surfaces that have the
same set of short curves one uses Minsky's product region theorem
\cite{minskyproduct}. To apply this, in addition to
Theorems~\ref{shortong} and~\ref{shortonL}, we need to estimate the
\Teich distance between the hyperbolic thick components of $\G_t$ and
$\L_t$, and also the difference between the Fenchel-Nielsen twist
coordinates corresponding to the short curves in the two surfaces.  In
Theorem~\ref{thm:markingsclose} and Corollary~\ref{equiv}, we show
that the \Teich distance between the corresponding thick components is
bounded.  In Theorem~\ref{thm:twistonL}, we estimate the twist of
$\np$ and $\nm$ around $\a$ at $\L_t$. Combined with the analogous
estimate for $\G_t$ proved in \cite{rafi2}, we are able to show that
the contribution to the \Teich distance between $\G_t$ and $\L_t$ from
the twisting is dominated by that from the lengths, leading to
\begin{introthm}[Theorem~\ref{thm:finalcomparison}]
\label{distance}
The \Teich distance between $\G_t$ and $\L_t$ is given by
$$d_{\ts}(\G_t,\L_t) = \frac{1}{2}\log \max \Big\{
\frac{l_{\mathcal G_t}(\a)}{l_{\L_t}(\a)} \Big\} \pm O(1),$$
\noindent where the maximum is taken over all curves $\a$ that are short in
$\mathcal G_t$.  
\end{introthm}

Theorems~\ref{shortong},~\ref{shortonL}, and~\ref{distance}
enable us to construct the various examples alluded to above.  Because
$K_t(\a)$ can become arbitrarily large while $D_t(\a)$ remains
bounded, it follows that $\G_t$ and $\L_t$ do not always remain a
bounded distance apart.  However, in the case in which $S$ is a
once-punctured torus or four-times-punctured sphere, it turns out that
the quantity $K_t(\a)$ is always bounded and therefore that the two
paths are always within bounded distance of each other. These ideas
are taken further in~\cite{crs}, where we show that $\L_t$ is a \Teich
quasi-geodesic.

\medskip

The greater part of the work of this paper is contained in the proof
of Theorem~\ref{shortonL}. It is carried out in several steps. First,
using the derivative formulae mentioned above, we show in
Theorem~\ref{thm:shortonL} that the length may be estimated by a
formula identical to that in Theorem~\ref{shortonL}, except that
$K_t(\a)$ is replaced by another geometric quantity
$$H_t(\a)=  
\sup_{\beta \in {\cal B}} \frac{l_{q_t}(\beta)}{l_{q_t}(\a)}. 
$$
Here $\cal B$ are those pants curves in a short pants decomposition
of $ {\L_t}$ (see Section~\ref{sec:shortmarkings}) which are boundaries
of pants adjacent to $\a$, while as above $ l_{q_t}$ denotes length in
the singular Euclidean metric associated to $q_t$.  
This is the content of Section~\ref{sec:shortlm}.

We now need to compare $H_t(\a)$ and $ K_t(\a)$.  From the definition
it is quite easy to show (Proposition~\ref{cor:shorteronL}) that
$H_t(\a) \succ K_t(\a)$.  In particular, it follows that a curve that
is short in $\G_t$ is at least as short in $\L_t$.  Next, we show in
Proposition~\ref{bddbelow} that on a subsurface whose injectivity
radius is bounded below in $\mathcal G_t$, the injectivity radius with
respect to $\L_t$ is also bounded below, perhaps by a smaller
constant.  The main point in the proof is
Proposition~\ref{lem:minimalonQ}, which shows that the hyperbolic
metric on $\L_t$ not only minimizes the sum of lengths $l({\nu^+_t}) +
l({\nu^-_t})$, but also, up to multiplicative error, it minimizes the
contribution of $l({\nu^+_t}) + l({\nu^-_t})$ to each thick component
 of the thick-thin
decomposition of $\G_t$. This proves Theorem~\ref{shortonboth}.

Having set up a one-to-one correspondence between the thick components
of $\G_t$ and $\L_t$, we show in Theorem~\ref{thm:markingsclose} and
Corollary~\ref{equiv} that the \Teich distance between corresponding
thick components is bounded.  Finally we are able to prove in
Proposition~\ref{lem:HandK} that $H_t(\a) \prec K_t(\a)$, completing
the proof of Theorem~\ref{shortonL}.

\medskip

Prior to this paper, the only results related to the relative behavior
of $\G$ and $\L$ were some partial results about their behavior at
infinity.  Results of Masur~\cite{masur} (for Teichm\"uller geodesics)
and of D\'{\i}az and Series~\cite{diazseries} (for lines of minima)
show that if either $\nu^{\pm}$ are supported on closed curves, or if
$\nu^{\pm}$ are uniquely ergodic, then $\G$ and $\L$ limit on the same
points in the Thurston boundary of $\teich(S)$.  In general, the
question of the behavior at infinity remains unresolved, but see
also~\cite{Lenzhen} which shows that there are Teichm\"uller geodesics
$\G$ which do not converge in Thurston's compactification of $\ts$. It
is not hard to apply the results of this paper to show the same is
true of lines of minima in Lenzhen's example; we hope to explore this
in more detail elsewhere.

\medskip

The motivation for our approach stems in part from a central
ingredient of the proof of the ending lamination theorem~\cite{elc}.
Suppose that $N$ is a hyperbolic $3$-manifold homeomorphic to $S
\times \RR$. The ending lamination theorem states that $N$ is
completely determined by the asymptotic invariants of its two ends.  A
key step is to show that if these end invariants are induced by the
laminations $\nu^+$ and $\nu^-$, then the curves on $S$ which have
short geodesic representatives in $N$ can be characterized in terms of
their combinatorial relationship to $\nu^+$ and $\nu^-$.  (The
relationship is expressed using the complex of curves of $S$, details
of which are not needed in what follows. Roughly speaking, a curve is
short in $N$ if and only if the distance between the projections of
$\nu^+$ and $\nu^-$ to some subsurface $Y \subset S$ is large in the
curve complex of $Y$.) In~\cite{rafi1}, Rafi found a similar
combinatorial characterization which shows that the curves which are
short in $N$ are almost, but not quite, the same as those curves which
become short along $\G(\np,\nm)$.  Our definition of $K_t$ is closely related
to Rafi's study~\cite{rafi3} of the relationship between the
thick-thin decomposition of a hyperbolic surface $S$ and a quadratic
differential metric on the same surface.  The relationship between
these two metrics plays a key role throughout the paper.

\medskip

The paper is organized as follows. In Section 2, we recall some
background facts about lines of minima and \Teich geodesics.  In
Section 3, we prove Theorem~\ref{thm:warmup} mentioned above, which
states that if both $\G_t$ and $\L_t$ are contained in the thick part
of $\ts$, then the \Teich distance between them is bounded. We hope
that treating this special case separately early on will give some
intuition about what needs to be done in general.  In Section 4, we
review twists and Fenchel-Nielsen coordinates and in Section 5, after
reviewing some fundamental facts about quadratic differential metrics,
we derive Theorem~\ref{shortong} and state the estimates for twists
about short curves proved in \cite{rafi2}. In Section 6, we prove
Theorem~\ref{thm:shortonL} and derive estimates for twists about the
short curves. Finally, in Section 7, we prove Theorems~
\ref{shortonboth}, \ref{shortonL} and ~\ref{distance}. Throughout the
paper, we make use of several basic length estimates on hyperbolic
surfaces.  The proofs, being somewhat long but relatively
straightforward, are relegated to the Appendix.

\subsection{Acknowledgments} We thank the referee for carefully reading
the manuscript and leading us to clarify the exposition. 

\section{Preliminaries}
\label{sec:prelims}
Throughout, $S$ is an orientable hyperbolizable surface of finite
type, possibly with punctures but with no other boundary.
\subsection{Thick-thin decomposition}
\label{sec:ttdecomp}
Let $\S$ denote the set of free homotopy classes of non-peripheral,
non-trivial simple closed curves on $S$.  If $(S, \s)$ is a surface
with hyperbolic metric $\s$ and $\a \in \S $, we write $l_\s(\a)$ for
the hyperbolic length of the unique geodesic representative of $\a$
with respect to $\s$.  The Margulis lemma provides a universal
constant $\ep_{\M}>0$ such that all components of the $\ep_{\M}$-thin
part of $(S,\s)$ (i.e.~the subset of $S$ where the injectivity radius
is less than $\ep_{\M}$) are horocyclic neighborhoods of cusps or
annular collars about short geodesics. The $\ep_{\M}$-thick part of
the surface is the complement of the thin part.

For our purposes, it is necessary to choose a constant $\ep_0>0$
sufficiently smaller than $\ep_{\M}$, in order that the
$\ep_0$-thick-thin decomposition of a surface satisfies certain
geometric conditions.  These conditions will be mentioned when the
context arises, but we assume that $\ep_0$ has been chosen once and
for all so that these conditions are met. If $l_\s(\a) < \ep_0$, we
shall say that $\a$ is {\em extremely short} in $\s$.
\subsection{Notation}
Since we will be dealing mainly with coarse estimates, we want to
avoid keeping track of constants which are universal, in that they do
not depend on any specific metric or curve under discussion.  For
functions $f,g$ we write $f\asymp g$ and $f \mul g$ to mean
respectively, that there are constants $c>1,C>0$, depending only on
the topology of $S$ and the fixed constant $\ep_0$, such that
$$\frac{1}{c} g(x) - C \leq f(x) \leq c g(x) +C \ \mbox{ and } \ 
\frac{1}{c} g(x) \leq f(x) \leq c g(x) .$$
The symbols ${\prec}$,
$\stackrel{{}_\ast}{\prec}$, ${\succ}$ , $\stackrel{{}_\ast}{\succ}$
are defined similarly.  For a positive quantity $X$, we often write
$X=O(1)$ instead of $X \prec 1$ to indicate $X$ is bounded above by a
constant depending only on the topology of $S$ and $\ep_0$, and more
generally we write $X = O(Y)$ to mean that $X/Y = O(1)$ for a positive
function $Y$.

\subsection{Measured laminations}
We denote the space of measured laminations on $S$ by $\ML$.  Given
any hyperbolic metric $\s$ on $S$, a measured lamination $\xi \in \ML$
can be realized as a geodesic measured lamination with respect to
$\s$. The hyperbolic length function extends by linearity and
continuity to $\ML$; we write $l_\s(\xi)$ for the hyperbolic length of
a lamination $\xi \in \ML$.  The geometric intersection number
$i(\a,\b)$ of curves $\a,\beta \in \S$ also extends continuously to
$\ML$.  Laminations $\mu,\nu \in \ML$ are said to \emph{fill up} $S$
if $i(\mu,\xi) + i(\nu,\xi) >0$ for all $\xi \in \ML$. For $\xi \in
\ML$, we denote the underlying leaves by $|\xi|$.

\subsection{Teichm\"uller space}
\label{sec:teichspace}
The Teichm\"uller space $\teich(S)$ of $S$ is the space of all
conformal structures on $S$ up to isotopy. The Teichm\"uller distance
$d_{\ts} (\Sigma,\Sigma')$ between two marked Riemann surfaces
$\Sigma, \Sigma' \in \ts$ is $[\log K]/2$, where $K$ is the smallest
quasiconformal constant of a homeomorphism from $\Sigma$ to $\Sigma'$
which is isotopic to the identity.

Each conformal structure $\Sigma \in \teich(S)$ is uniformized by a
unique hyperbolic structure $\s$, and conversely, each hyperbolic
structure $\s$ has an underlying conformal structure $\Sigma$. Thus,
we also consider $\ts$ to be the space of all hyperbolic metrics on
$S$ up to isotopy. The \emph{thick part} $\tthick$ of $\teich(S)$ will
be defined as the subset of all hyperbolic metrics such that every
closed geodesic has length bounded below by the constant $\ep_0$.

\subsection{Kerckhoff lines of minima}
Suppose that $\np,\nm \in \ML$ fill up $S$.  Kerckhoff~\cite{lminima}
showed that the length function 
$$\sigma \mapsto l_\s(\nu^{+}) +l_\s(\nu^{-})$$
has a global minimum
on $\teich(S)$ at a unique surface $\L_0$.  Moreover, as $t$ varies in
$(-\infty,\infty)$, the surface $\L_t \in \teich(S)$ that realizes the global minimum of
$l(\nu_t^+)+l(\nu_t^-)$ for the weighted laminations $\nu^{+}_t=
e^{t}\np$ and $\nu^{-}_t= e^{-t}\nm$ varies continuously with $t$ and
traces out a path $t \mapsto \L_t$ called the {\em line of minima
  $\L(\np,\nm)$ of }$\nu^{\pm}$.

\subsection{Quadratic differentials}
\label{sec:quaddiffs}
We give a brief summary of facts about quadratic differentials that we
use and refer the reader to ~\cite{gardiner, strebel} for a detailed
and comprehensive background. Let $\Sigma$ be a Riemann surface and
$q$ a quadratic differential on $\Sigma$ which is holomorphic, except
possibly at punctures, where $q$ may have a pole of order one. This
ensures that the area of $\Sigma$ with respect to the area element
$|q(z) dz^2|$ is finite, and we normalize so that the area is $1$. Let
$\cal Q(\Sigma)$ be the space of all such meromorphic quadratic
differentials on $\Sigma$.

The zeros and poles of $q$ are called {\em critical points}. Away from
the critical points, we have two mutually orthogonal line fields
defined respectively by the conditions that $\mbox{Im} [\sqrt{q(z)}
dz]$ is zero and $\mbox{Re}[\sqrt{q(z)}dz]$ is zero.  This defines a
pair of measured singular foliations on $\Sigma$ with singularities at
the critical points of $q$, respectively called the {\em horizontal
  foliation} $\cal H_q$ and the {\em vertical foliation} $\cal V_q$.
The measures on these foliations are determined by integrating the
line element $|\sqrt{q(z)} dz|$. More precisely, for a curve $\eta$,
its horizontal and vertical measures are given respectively by
$$h_q(\eta) = \int_\eta | \mbox{Re} [\sqrt{q(z)} dz]|, \ v_q(\eta) =
\int_\eta | \mbox{Im} [\sqrt{q(z)} dz]|.$$
We call $h_q(\eta)$ and
$v_q(\eta)$ respectively, the {\em horizontal length} and the {\em
  vertical length} of $\eta$.

Every essential simple closed curve $\gamma$ in $(S,q)$ has a unique
$q$-geodesic representative, unless it is in a family of closed
Euclidean geodesics foliating an annulus whose interior contains no
singularities. We denote the $q$-geodesic length of $\gamma$ by
$l_q(\gamma)$. It satisfies the following inequalities:
\begin{equation}
\label{handv}
[h_q(\gamma) + v_q(\gamma)]/\sqrt{2} \leq l_q(\gamma) \leq
h_q(\gamma) + v_q(\gamma).
\end{equation}
By definition of intersection numbers for measured foliations, we have
$v_q(\gamma) = i(\cal H_q,\gamma)$ and $h_q(\gamma) = i(\cal
V_q,\gamma)$,
so Equation~(\ref{handv}) implies
\begin{equation}
\label{eqn:qlengthapprox}
l_q(\gamma) \mul i(\cal V_q,\gamma) + i(\cal H_q,\gamma).
 \end{equation}
This approximation will be used repeatedly.
\subsection{\Teich geodesics}

Suppose that $\Sigma,\Sigma' \in \teich(S)$ are marked Riemann
surfaces with $d_{\T(S)}(\Sigma,\Sigma') = d$. Then there is a unique
quadratic differential $q$ on $\Sigma$ such that the conformal
structure on $\Sigma'$ is obtained from that of $\Sigma$ by expanding
in the horizontal direction of $q$ by a factor $e^{d}$ and contracting
in the vertical direction by $e^{-d}$.  The homeomorphism which
realizes this is called the Teichm\"uller map from $\Sigma$ to
$\Sigma'$, and has quasiconformal distortion $e^{2d}$. The
$1$-parameter family of quadratic differentials $q_t$ whose horizontal
and vertical foliations are respectively, $e^t \cal H_q$ and $e^{-t}
\cal V_q$, for $0 \leq t \leq d$, define the geodesic path from
$\Sigma$ to $\Sigma'$ with respect to the \Teich metric.

Gardiner and Masur~\cite{G-M} showed that for any pair of measured
laminations $ \nu^+, \nu^- \in \ML$ which fill up $S$ and such that $i(\nu^+,
\nu^-)=1$, there is a unique Riemann surface $\Sigma \in \ts$ and a
unique quadratic differential $q \in \Q(\Sigma)$ whose horizontal and
vertical foliations are $ \nu^+, \nu^-$ respectively. (This uses the
one-to-one correspondence between the space of measured laminations
and the space of measured foliations.)  For $t \in \RR$ set
$$\nu^+_t = e^{t}\nu^+, \quad \nu^-_t= e^{-t} \nu^-$$ 
and let $\G_t$ and $q_t$ be the corresponding Riemann surface
and quadratic differential. The path $t \mapsto \G_t$ defines a
Teichm\"uller geodesic which we denote $\G=\G( \nu^+, \nu^-)$.  
We abuse notation and use $\G_t$ to also denote the hyperbolic metric
that uniformizes the Riemann surface $\G_t$.

 \subsection{The balance time}
 \label{sec:balancetime}
 Let $\a \in \S$.  We say $\a$ is \emph{vertical} along
 $\G(\np,\nm)$ if its intersection $i(\a,\nu^-)$ with the vertical
 foliation $\nm$ vanishes. In this case, $\a$ can be realized as a
 union of leaves of the vertical foliation. Similarly, $\a$ is
 \emph{horizontal} if $i(\a,\nu^+) = 0$.  Mostly we shall be dealing
 with curves $\a$ which are neither horizontal nor vertical.  In this
 case, there is always a unique time $t_\a$ at which
 $i(\a,\nu^+_{t_\a}) = i(\a,\nu^-_{t_\a})$. We call $ t_\a$ the
 \emph{balance time} of $\a$. The length of $\a$ with respect to
 $\G_t$ is approximately convex along $\G$ and is close to its minimum
 at $ t_\a$, see~\cite{rafi2} Theorem 3.1.  Our estimation of the
 hyperbolic lengths of short curves will mainly be made relative to
 their balance time.

 \section{Comparison on the thick part of Teichm\"uller space}
\label{sec:thickparts}
In this section we prove Theorem~\ref{thm:warmup}, which states that
if $\G_t$ and $\L_t$ are in the $\ep_0$-thick part $\T_{thick}(S)$ of
\Teich space, then the \Teich distance between them is uniformly
bounded by a constant that depends only on the topology of $S$ and
$\ep_0$. The idea is to first approximate the length of a curve
$\zeta$ for any $\s \in \T_{thick}(S)$ by its intersection with what
we call a {\em short marking} for $\s$, and then to compare the short
markings for $\G_t$ and $\L_t$. This method will be extended in
Section~\ref{sec:finalcompare} when we consider the \Teich distance
between $\G_t$ and $\L_t$ in general. We begin with some definitions.
\subsection{Short markings}
\label{sec:shortmarkings}
We call a maximal collection of pairwise disjoint, homotopically
distinct, non-peripheral, non-trivial simple closed curves on $S$, a
\emph{pants curve system} on $S$. The terminology is due to the fact
that the complementary components are pairs of pants, i.e., three
holed spheres (in which some boundary components may be punctures).
Our notion of a {\em marking} is motivated by \cite{masurminskyII}:
\begin{definition}
\label{defn:markings}
A marking $M$ on a surface $S$ is a system of pants curves 
$\a_1,\ldots, \a_k$ and simple closed curves $\delta_{\a_1},\ldots,
\delta_{\a_k}$ such that
$$
\begin{cases}
  i(\a_i, \delta_{\a_j}) = 0 &\text{if \quad $i\neq j$} \\
  i(\a_i,\delta_{\a_i})=2 & \text{if \quad two distinct pairs
    of pants are adjacent along $\a_i$}\\
  i(\a_i,\delta_{\a_i})=1 &\text{if \quad $\a_i$ is adjacent to only a
    single pair of pants.}
\end{cases}
$$
We call $\delta_{\a_i}$ the dual curve of $\a_i$. 
\end{definition}
In the second case, $\a_i \cup \delta_{\a_i}$ fill a four-holed sphere
(that is, a regular neighborhood of $\a_i \cup \delta_{\a_i}$ is
homeomorphic to a four-holed sphere) and in the third case $\a_i \cup
\delta_{\a_i}$ fill a one-holed torus.  It is easy to see that any two
markings which have the same pants system $\cal P$ have dual curves
which differ only by twists and half-twists around the curves in $\cal
P$.

The following well-known lemma states that for any hyperbolic metric,
one can always choose a pants system whose length is universally
bounded:
\begin{lemma}
\label{lem:Bers}
(Bers~\cite{bers1}) There exists a constant $L>0$ such that for every
$\s \in \ts$ there is a pants curve system $\cal P$ with the property
that $l_{\s}(\a) < L$ for every $\a \in \cal P$.
\end{lemma}
\noindent If the boundary curves of a pair of pants have bounded
length as in Bers's lemma, the geometry of a pair of pants satisfies
the following (for a proof, see Appendix):
\begin{lemma}
\label{lem:width}
Let $P$ be a totally geodesic pair of pants with boundary curves
$\a_1,\a_2,\a_3$ of lengths $l(\a_i)<L$ for $i=1,2,3$.  Then the
common perpendicular of $\a_i,\a_j$ (where possibly $i=j$) has length
$$\log \frac1{l(\a_i)} + \log \frac 1{l(\a_j)} \pm O(1), $$ 
where the bound on the error depends only on $L$.
\end{lemma} 

We will say that a pants curve system as in Lemma~\ref{lem:Bers} is
{\em short} in $(S,\s)$.  A {\em short marking} $M_\s$ for $\s$ is a
short pants system together with a dual system chosen so that each
dual curve $\delta_{\a_i}$ is the shortest among all possible dual
curves. For a given pants curve, notice that there may be more than
one shortest dual curve, in which case any choice will suffice.  Also
notice that \emph{not all curves in a short marking are necessarily
  short}; if a pants curve is very short, then the corresponding dual
curve will be very long. More precisely, we have the following easy
consequence of Lemma~\ref{lem:width}:
 \begin{corollary}
\label{cor:duallength}
Let $M_\s$ be a short marking for $\s$ and let $\a, \delta_\a \in
M_\s$ be a pants curve and its dual. Then
$$l_\s(\delta_\a) = i(\delta_\a, \a) \cdot 2 \log \frac1{l_\s(\a)} \pm O(1).$$
\end{corollary} 
Observe that if $\s \in \tthick$, since the length of every curve in
$M_\s$ is uniformly bounded below, it follows from
Lemma~\ref{lem:Bers} and Corollary~\ref{cor:duallength} that the
length of every curve in $M_\s$ is also uniformly bounded above. Thus,
if $\s \in \tthick$, we have $l_\s(M_\s) \mul 1$.
 

For surfaces in the thick part of $\ts$, short markings coarsely
determine the geometry.  We express this in the following proposition
whose proof can be found in the proof of Lemma 4.7 in
~\cite{minskytop} (see also the proof of
Proposition~\ref{prop:capmarking1} below).  If $M$ is a marking and
$\xi \in \ML$, we write
$$i(M, \xi) = \sum_{\g \in M} i(\g,\xi).$$
\begin{proposition}
\label{prop:capmarking0}
Let $M_\s$ be a short marking for $\s \in \tthick$.  Then for 
any $\zeta \in {\cal S}$,
$$
l_\s (\zeta) \mul i(M_\s, \zeta).$$
\end{proposition}
\noindent
Since both length and intersection number scale linearly with weights
of simple closed curves, it follows that:
\begin{proposition}
\label{prop:capmarking}
Let $M_\s$ be a short marking for $\s \in \tthick$.  Then for any
$\xi \in \ML$,
$$l_\s (\xi) \mul i(M_\s, \xi).$$
\end{proposition}
\subsection{Comparison on the thick part}
We use the estimate in Proposition~\ref{prop:capmarking} to compare
$\G_t$ and $\L_t$ in the thick part of $\ts$.  The following
well-known lemma is proved in greater generality in \cite{rafi3} (see
Theorem~\ref{thm:scalefactor}(ii) below).  Recall that $\Sigma$
denotes the conformal structure associated to the metric $\s$.
\begin{lemma} 
\label{lemma:compare} Suppose that $\s \in \tthick$  and $q \in
{\cal Q}(\Sigma)$. Then for every $\zeta \in \S$,
$$l_{\s}(\zeta) \mul l_{q}(\zeta).$$ 
\end{lemma}

\begin{theorem}
\label{thm:warmup}
If $\G_t, \L_t \in \tthick$ then $ d_{\ts}(\G_t,\L_t) = O(1).$
\end{theorem}
\begin{proof}
  By Lemma~\ref{lemma:compare}, Equation~(\ref{eqn:qlengthapprox}),
  and Proposition~\ref{prop:capmarking}, we have
\begin{align}
\label{eqn:forg} 
l_{\mathcal G_t}(M_{\mathcal G_t}) &\mul l_{q_t}(M_{\mathcal G_t})
\mul i(M_{\mathcal G_t},\nu^+_t) +
i(M_{\mathcal G_t},\nu^-_t) \\
&\mul l_{\mathcal G_t}(\nu^+_t) + l_{\mathcal G_t}(\nu^-_t). \nonumber
\end{align}
Now, since $\L_t$ minimizes $l_\s(\nu^+_t) + l_\s(\nu^-_t)$ over all
$\s \in \ts$, we have
$$
l_{\mathcal G_t}(\nu^+_t) + l_{\mathcal G_t}(\nu^-_t) \geq  
l_{\L_t}(\nu^+_t) + l_{\L_t}(\nu^-_t).
$$
Reversing the sequence of estimates in Equation~(\ref{eqn:forg}),
we get
\begin{align*}
l_{\L_t}(\nu^+_t) + l_{\L_t}(\nu^-_t)  
 &\mul i(M_{\L_t},\nu^+_t) + i(M_{\L_t},\nu^-_t) \mul l_{q_t}(M_{\L_t}) \\ 
 &\mul l_{\mathcal G_t}(M_{\L_t}). 
\end{align*}
Putting together the preceding three equations, we have
$$l_{\mathcal G_t}(M_{\mathcal G_t}) \stackrel{{}_\ast}{\succ}
l_{\mathcal G_t}(M_{\L_t}).$$
Since $\mathcal G_t \in \tthick$, it
follows from the observation following Corollary~\ref{cor:duallength}
that
\begin{equation}
\label{eqn:n}
l_{\mathcal G_t}(M_{\L_t}) \stackrel{{}_\ast}{\prec} 1.
\end{equation}
Notice also that $l_{\L_t}(M_{\L_t}) \stackrel{{}_\ast}{\prec} 1$.
Lemma 4.7 of~\cite{minskytop} implies that for any given $B>0$, the
diameter of the set $\{\, \s \in \ts : l_\s(M_{\L_t}) < B \, \}$, with
respect to the \Teich distance, is bounded above by a constant that
depends only on $B$. Thus it follows that $d_{\ts}(\mathcal G_t,\L_t)
= O(1).$
\end{proof}

\section{Twists and Fenchel-Nielsen  coordinates}
\label{sec:twistsFN}
 
In order to compare surfaces in the thin part of Teichm\"uller space,
our main tool will be Minsky's product region
theorem~\cite{minskyproduct}. This uses Fenchel-Nielsen coordinates to
give a nice coarse expression for Teichm\"uller distance between
surfaces which have common thin parts.  To state the results
precisely, we first discuss twists and Fenchel-Nielsen coordinates.

\subsection{Twists in hyperbolic metrics}
\label{sec:twistshyperbolic}
\begin{figure}[htb]
\centerline{\epsfbox{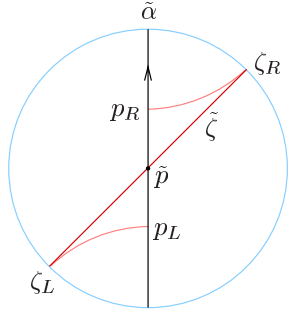}}
\caption{Defining the twist of $\zeta$ around $\a$.}
 \label{fig:twist}
\end{figure}

There are various ways to define the twist of one curve around
another, all of which differ by factors unimportant to us here. We
shall follow Minsky \cite{minskyproduct}.  Let $\s \in \ts$ be a
hyperbolic metric and let $\a $ be an oriented simple closed geodesic
on $(S,\s)$. Let $\zeta$ be a simple geodesic that intersects $\a$
transversely and let $p$ be a point of intersection. In the universal
cover $\HH^2$, a lift $\tilde \zeta$ of $\zeta$ intersects a lift
$\tilde \a$ of $\a$ at a lift $\tilde p$ of $p$, and has endpoints
$\zeta_R, \zeta_L$ on $\dd_\infty \HH^2$ to the right and left of
$\tilde \a$, respectively (see Figure~\ref{fig:twist}).
Let $p_R, p_L$ be the orthogonal projections of $\zeta_R,\zeta_L$ to
$\tilde \a$ respectively. Then the twist of $\zeta$ around $\a$ at $p$
is defined as
$$
tw_\s(\zeta,\a,p) = \pm \frac{d_{\HH^2}(p_R,p_L)}{l_\s(\a)},$$
where the sign is $(+)$ if the direction from $p_L$ to $p_R$ coincides
with the orientation of $\tilde \a$ and $(-)$ if it is opposite.  For
any other point $q \in \zeta \cap \a$, the twist satisfies
(\cite[Lemma 3.1]{minskyproduct})
$$|tw_\s(\zeta,\a, q)-tw_\s(\zeta,\a,p)| \leq 1.$$
To obtain a number that is independent of the point of intersection, 
Minsky defines 
$$tw_\s(\zeta,\a) = \min_{p \in \zeta \cap \a} tw_\s(\zeta,\a,p).$$
For convenience, we write $Tw_\s(\zeta,\a)$ for $|tw_\s(\zeta,\a)|$.

Note that the definition of twist is valid even if the simple geodesic
$\zeta$ is not closed, because the inequality 
$$|tw_\s(\zeta,\a, q)-tw_\s(\zeta,\a,p)| \leq 1$$
depends only on the
fact that different lifts of $\zeta$ are disjoint.  Thus if $\nu$ is a
measured geodesic lamination that intersects $\a$ transversely, we can
define the twist of $\nu$ around $\a$ by taking the infimum of twists
over all leaves of $\nu$ that intersect $\a$.  We remark that although
we will be working with measured geodesic laminations, when defining
the twist, the measure is irrelevant, in other words the twist
$tw_\s(\nu,\a)$ depends only on the underlying lamination $| \nu|$.

\subsection{Fenchel-Nielsen  coordinates}
\label{sec:fncoords}
We define the Fenchel-Nielsen coordinates $$(l_{\s}(\a_i),
s_{\a_i}(\s))_{i=1}^k$$
associated to a pants curves system
$\a_1,\ldots,\a_k$ in the following standard way, see for example
\cite{minskyproduct}.  Suppose that $P$ is a pair of pants that is
the closure of a component of $S\setminus \{\a_1,\ldots,\a_k\}$.
By a {\em seam} of $P$, we mean a common perpendicular between two
distinct boundary components of $P$.  (Notice that the definition of
seam refers to the internal geometry of $P$ alone; two distinct
boundary curves of $P$ may project to the same curve on $S$.) Each
boundary curve of $P$ is bisected by the two points at which it meets
the two seams intersecting it.  We
first construct a base surface $\sigma_0 = \sigma_0 (l_{1}^0,\ldots,
l_{k}^0) $ in which the pants curve $\a_i$ has some specific choice of
length $l_{i}^0$.  Each $\a_i$ is adjacent to two (possibly
coincident) pairs of pants; we glue these two pants together in such a
way that seams incident on $\a_i$ from the two sides match up.  Since
the seams meet the pants curves orthogonally, they glue up to form a
collection of closed geodesics $\gamma_j$.  Any other structure
$\sigma \in \ts$ comes endowed with an associated homeomorphism
$h:\sigma_0 \to \sigma$.  The length coordinates of $\sigma$ are
defined by $ l_{\s}(\a_i)$. Let $\sigma_0( l_{1},\ldots,l_k)$ denote
the surface in which $\a_i$ (more precisely $h(\alpha_i)$) has length
$l_{i} = l_{\s}(\a_i)$, while the curves formed by gluing the new
seams are exactly the images $h(\gamma_j)$.  Now define
$\tau_{\a_i}(\s)$ to be the signed distance that one has to twist
around $\a_i$ to obtain $\s$ starting from $\sigma_0(
l_{1},\ldots,l_k)$, where the sign is determined relative to a fixed
orientation on $\sigma_0$ and hence on $\s$.  Finally we define the
twist coordinates of $\s$ by
$$s_{\a_i}(\s) = \frac{ {\tau}_{\a_i}(\s)}{l_{\s}(\a_i)} \in \RR.$$

\begin{lemma} (Minsky \cite{minskyproduct} Lemma 3.5) 
\label{lem:minskytwist}
For any lamination $\nu \in \ML$ that intersects $\a=\a_i$ and any two
metrics $\s, \s' \in \ts$,
$$ | (tw_\s(\nu,\a) - tw_{\s'}(\nu,\a)) -
(s_\a (\s)  -
 s_{\a} (\s'))| \leq 4.$$ 
\end{lemma}

In~\cite{minskyproduct}, the statement is only given for closed
curves. However the argument extends without change to laminations.
This is because the proof in \cite{minskyproduct} is based on the
observation that for any two simple closed curves $\zeta_1, \zeta_2$
intersecting $\a$, the difference $tw_\s(\zeta_1,\a)-
tw_\s(\zeta_2,\a)$ is a topological quantity, independent of $\s$, up
to a bounded error of $1$. More precisely, it follows from the proof
in \cite{minskyproduct} that if $\tilde S$ is the annular cover of $S$
corresponding to $\a$,  and if $\tilde \zeta_1$ and $\tilde \zeta_2$ are
respectively, lifts of $\zeta_1$ and $\zeta_2$ intersecting the core
$\tilde \a$ of $\tilde S$, then the difference $tw_\s(\zeta_1,\a)-
tw_\s(\zeta_2,\a)$ is the signed intersection of $\tilde \zeta_1$ and
$\tilde \zeta_2$ in the annulus $\tilde S$, up to a bounded error.
This topological characterization holds even when $\zeta_1$ and
$\zeta_2$ are simple geodesics which are not
necessarily closed.
 
  \subsection{Relative twist in an annulus}
  \label{sec:reltwist}

The above topological observation allows us to define the following: 
\begin{definition}
\label{defn:algint}
For any two laminations $\nu_1,\nu_2$ that intersect a curve $\a$,
define their algebraic intersection around $\a$ to be
$$ i_\a(\nu_1,\nu_2) = \inf_\s [tw_\s(\nu_1,\a) -
tw_\s(\nu_2,\a)],$$
where the infimum is taken over all possible surfaces $\s \in \ts$.
\end{definition}
\noindent Often we need only the absolute value:  
\begin{definition}
For any two laminations $\nu_1,\nu_2$ that intersect a curve $\a$,
define the relative twisting of $\nu_1,\nu_2$ around $\a$ to be 
$$d_\a(\nu_1,\nu_2) = |\,i_\a(\nu_1,\nu_2)\,|.$$
\end{definition}
\noindent Thus for any $\s \in \ts$, we have 
$$|tw_\s(\nu_1,\a) -
tw_\s(\nu_2,\a)|= d_\a(\nu_1,\nu_2) + O(1).$$  Notice that
$i_\a(\nu_1,\nu_2)$ and $d_\a(\nu_1,\nu_2) $ are independent of the
measures on $\nu_1,\nu_2$, depending only on the underlying
laminations $|\nu_1|$ and $|\nu_2|$.

Using Definition~\ref{defn:algint}, one sees easily that
\begin{equation}
\label{eqn:reltwist} 
i_\a(\nu_1,\nu_2) = i_\a(\nu_1,\xi) - i_\a(\nu_2,\xi) + O(1)
\end{equation}
for any curve $\xi$ transverse to $\a$. It is also easily seen that $
d_\a(\nu_1,\nu_2) $ agrees up to $O(1)$ with the definition of
subsurface distance between the projections of $|\nu_1|$ and $|\nu_2|$
to the annular cover of $S$ with core $\a$, as defined
in~\cite{masurminskyII} Section 2.4 and used throughout~\cite{rafi1,
  rafi2}.

\medskip

Another essentially equivalent way of measuring twist is to look at
the intersection with the shortest curve transverse to $\a$: 
\begin{lemma}
\label{lemma:twistintersect}
Let $\a$ be a pants curve and let $\delta_\a$ be a shortest dual
curve of $\a$ in some marking for $\s$. Then for any simple closed
curve $\zeta$ intersecting $\a$, we have $|tw_\s(\zeta,\a) -
i_{\a}(\zeta,\delta_\a) | = O(1)$.
\end{lemma}
\begin{proof}
Since $i_\a(\zeta ,\delta_\a) = tw_\s(\zeta,\a) - tw_\s(\delta_\a,
\a)$, up to a bounded error of $1$, it is sufficient to show that
$|tw_\s(\delta_\a,\a)| = O(1)$.

Let $\tilde \delta_\a$ be a lift of $\delta_\a$ in the universal cover
$\HH^2$ and let $\tilde \a, \tilde \a'$ be the two lifts of $\a$
containing the endpoints of $\tilde \delta_\a$ (see
Figure~\ref{fig:deltatwist}). Let $\eta$ be the perpendicular between
$\tilde \a$ and $\tilde \a'$ and let
\begin{figure}[htb]
\centerline{\epsfbox{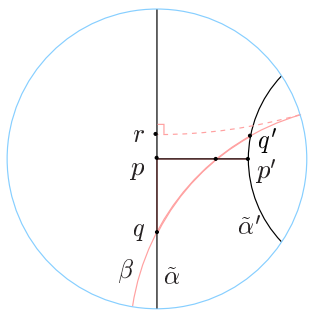}}
\caption{Bound on the twist of $\delta_\a$ around $\a$.}
 \label{fig:deltatwist}
\end{figure}
$p,p'$ be the endpoints of $\eta$ on $\tilde \a,\tilde \a'$,
respectively. Since $\delta_\a$ is the shortest dual curve, the
endpoints of $\tilde \delta_\a$ must be within distance $l_\s(\a)$
from $p,p'$.  Let $q,q'$ be points on $\tilde \a,\tilde \a'$ at
distance $l_\s(\a)$ from $p,p'$, respectively, on opposite sides of
$\eta$. It is easy to see that $|tw_\s(\delta_\a,\a)| \leq
|tw_\s(\beta,\a)|$, where $\beta$ is the geodesic through $q,q'$. Let
$r$ be the foot of the perpendicular as shown. 
As in the first part of the proof of Lemma 3.5
in~\cite{minskyproduct}, we note that since the images of $\tilde
\alpha'$ under the translation along $\tilde \alpha$ are disjoint, the
projection of $\tilde \alpha'$ on $\tilde \alpha$ has length at most
$l_\s(\a)$. Thus $l(pr) < l(pq) = l(\a)$. Hence, 
\begin{equation*}
|tw_\s(\beta,\a)|= 2 \, \frac{l(pq) + l(pr)}{l(\a)} < 4. \qedhere
\end{equation*}
\end{proof}

\subsection{The product region theorem}
\label{sec:productregion}
Let ${\cal A} \subset \S$ be a collection of disjoint, homotopically
distinct, simple closed curves on $S$ and let $\T_{thin}({\cal
  A},\ep_0) \subset \ts$ be the subset in which all curves $\a \in
{\cal A}$ have hyperbolic length at most $\ep_0$.  Extend ${\cal A}$
to a pants decomposition and define Fenchel-Nielsen coordinates as
above.  Let $S_{\cal A}$ denote the surface obtained from $S$ by
removing all the curves in $\cal A$ and replacing the resulting
boundary components by punctures.

Following~\cite{minskyproduct}, we now define a projection
$$
\Pi \colon \ts \to \T(S_{\cal A}) 
\times \HH_{\a_1} \times \ldots \times \HH_{\a_r},
$$
where $ {\cal A} = \{ \a_1 \ldots, \a_r\}$ and $\HH_{\a_i}$ is the
upper half-plane. The first component $\Pi_0$ which maps to $
\T(S_{\cal A})$ is defined by forgetting the coordinates of the curves
in $\cal A$ and keeping the same Fenchel-Nielsen coordinates for the
remaining surface.  For $\a \in {\cal A}$ define $\Pi_{\a} \colon \ts
\to \HH_{\a}$ by
$$\Pi_{\a}(\s) = s_{\a}(\s) + i/ l_{\s}(\a) \in \HH_{\a}.$$
Let
$d_{\HH_\a}$ be half the usual hyperbolic metric on $\HH_\a$.
Minsky's product region theorem states that, up to bounded additive
error, Teichm\"uller distance on $\T_{ thin}({\cal A},\ep_0) $ is
equal to the sup metric on
$$
\T(S_{\cal A}) \times \HH_{\a_1} \times \ldots \times \HH_{\a_r}.
$$

\begin{theorem} [Minsky \cite{minskyproduct}]
 \label{thm:minskyproduct} Let $\s,\t \in \T_{thin}({\cal
A},\ep_0) $.
Then $$d_{\ts}(\s,\t) =  \max_{\a \in \cal A} \{ d_{\T(S_{\cal
A})}(\Pi_0(\s),\Pi_0(\t)), d_{\HH_{\a}}(\Pi_{\a}(\s), \Pi_{\a}(\t) )
\} \pm O(1).$$
\end{theorem}
\noindent We remark that Minsky makes several assumptions on the size of
$\ep_0$ in order to prove the above theorem. We may assume that our
initial choice of $\ep_0$ satisfies these assumptions. Recall from
Section~\ref{sec:ttdecomp} that a curve $\a \in \S$ is said to be
extremely short if $l_{\sigma}(\alpha) < \ep_0$.  The distance between
the projections to $\HH_\a$ can be approximated as follows:
\begin{lemma}
 \label{lemma:hypdist} Suppose $\a \in {\cal S}$ is extremely short in both
 $\s,\t$. Then $$
 \exp2{d_{\HH_{\a}}(\Pi_{\a}(\s),\Pi_{\a}(\t))} \mul
 \max \left\{l_\s l_\t |s_{\a}(\s) - s_{\a}(\t)|^2, l_\t/l_\s,l_\s/l_\t \right\}
$$
where  $l_{\s}= l_\s(\a)$ and $l_\t= l_\t(\a)$.
\end{lemma}
\begin{proof} This is a simple calculation using the formula 
$$  \cosh  2d_{\HH}(z_1,z_2)  = 1 + \frac{|z_1-z_2|^2}{2 \,
\mbox{Im} \,z_1 \mbox{Im}\, z_2}   ,$$
where $d_{\HH}$ is as above half the usual hyperbolic distance in $\HH$.
\end{proof}

The lemma reveals a useful fact from hyperbolic geometry: unless the
difference $|x-x'|$ is extremely large, the distance between two
points $x+iy,x'+iy' \in \HH$ is dominated by $|\log [y/y']|$.  In our
situation this means that unless the difference between the twist
coordinates is extremely large in comparison to the lengths, their
contribution to hyperbolic distance can be neglected.  We quantify
this in the following useful corollary which shows that as long as
each twist coordinate is bounded by $O(1/l)$, it can be neglected when
estimating the contribution to \Teich distance coming from the same
short curve in two surfaces.

\begin{corollary}
\label{cor:hypcomparison} Suppose that $\s_1,\s_2 \in
\T_{thin}(\a,\ep_0)$ 
and that,
for $i=1,2$, $\nu \in \ML$ satisfies
$$Tw_{\s_i}(\nu,\a)\, l_{\s_i}(\a) = O(1).$$ 
Then
$$
d_{\HH_{\a}}(\Pi_{\a}(\s_1),\Pi_{\a}(\s_2) ) = \frac{1}{2}
\left| \log \frac{ l_{\s_1}(\a)}{l_{\s_2}(\a)} \right| \pm O(1).
$$
\end{corollary}
\begin{proof} This follows easily from Lemma~\ref{lemma:hypdist}, using
Lemma~\ref{lem:minskytwist} to approximate $s_{\a}(\s_1)-s_{\a}(\s_2)$ by 
$tw_{\s_1}(\nu,\a) - tw_{\s_2}(\nu,\a)$.  (Cross terms such as 
$s_\a(\s_1) l_{\s_2}(\a)$ may be rearranged as 
$[s_\a(\s_1) l_{\s_1}(\a)][l_{\s_2}(\a)/l_{\s_1}(\a)]$.) Note that
the multiplicative error in Lemma~\ref{lemma:hypdist} translates to
an additive error in the distance.
\end{proof}
\section{Short curves along Teichm\"uller geodesics}
\label{sec:shortonG}
In this section we prove Theorem~\ref{shortong}, stated more precisely
as Theorem~\ref{prop:shortingeo}. This gives a combinatorial estimate
for the hyperbolic length of an extremely short curve along the \Teich
geodesic $\G(\np,\nm)$.  We also recall the estimate for the twist of
$\nu^\pm$ around such curves (Theorem~\ref{prop:quadtwistest}) proved
in \cite{rafi2}.

Deriving the length estimate is largely a matter of putting together
results proved in \cite{minskyharmonic, rafi1, rafi3}. Both estimates
are made by a careful study of the relationship between the hyperbolic
metric $\mathcal G_t$ and the quadratic differential metric $q_t$.  As
indicated in the Introduction, there are two distinct reasons why a
curve $\a$ may become extremely short: one is that the relative
twisting of $\nu^+$ and $\nu^-$ around $\a$ may be very large, the
other that $\nu^+$ and $\nu^-$ may have large relative complexity in
$S \setminus \a$.  We express the latter in terms of a scale factor
which controls the relationship between the quadratic differential and
hyperbolic metrics on the components of the thick part of $\mathcal
G_t$ adjacent to $\a$, as made precise in Rafi's thick-thin
decomposition for quadratic differentials,
Theorem~\ref{thm:scalefactor}.

We begin in Sections~\ref{sec:annuli} and \ref{sec:modofann} with some essential
facts about the geometry of annuli with respect to quadratic
differential metrics.

\subsection{Annuli in quadratic differential metrics}
\label{sec:annuli}
We review the notions of {\em flat} and {\em expanding} annuli from
\cite{minskyproduct}. These concepts provide the framework with which
to analyze short curves.

Let $\Sigma \in \ts$ and let $q$ be a quadratic differential on
$\Sigma$.  Let $\gamma$ be a piecewise smooth curve in $(S,q)$.  At a
smooth point $p$, the curvature $\kappa(p)$ is well-defined, up to a
choice of sign. If $\gamma$ is the boundary component of a subsurface
$Y$, we choose the sign to be positive if the acceleration vector at
$p$ points into $Y$. At a singular point $P$, although the curvature
is not defined, we shall say $\gamma$ is non-negatively curved at $P$
if the interior angle $\theta(P)$ is at most $\pi$ and say it is
non-positively curved at $P$ if $\theta(P)$ is at least $\pi$.  By
interior angle, we mean the angle that is on the same side of $\gamma$
as $Y$. We say $\gamma$ is monotonically curved with respect to $Y$
either if the curvature is non-negative at every point, or
non-positive at every point. The total curvature of $\gamma$ is given
by
$$\kappa_Y(\gamma) = \int_{\gamma} \kappa(p) + \sum[\pi - \theta(P)],$$
where the sum is taken over all singular points $P$ on $\gamma$.
The Gauss-Bonnet theorem gives
\begin{equation}
\label{eqn:gaussbonnet}
 \sum \kappa_Y(\gamma) - \pi \sum \mbox{ord}P= 
2 \pi \chi(Y), 
\end{equation}
where $\mbox{ord}P$ is the order of the zero at $P$, the first sum is
over all boundary components $\gamma$ of $Y$, and the second sum is
over the zeros $P$ of $q$ in the interior of $Y$.

Let $A$ be an annulus in $(S,q)$ with piecewise smooth boundary. The
following definitions are due to Minsky \cite{minskyharmonic}. We say
$A$ is {\em regular} if both boundary components $\dd_0, \dd_1$ are
monotonically curved with respect to $A$ and if $\dd_0$, $\dd_1$ are
$q$-equidistant from each other. Suppose that $A$ is a regular annulus
such that $\kappa_A(\dd_0) \leq 0$. We say $A$ is an {\em expanding
  annulus} if $\kappa_A(\dd_0)<0$ and we call $\dd_0$ the {\em inner}
boundary and $\dd_1$ the {\em outer} boundary.  Expanding annuli are
exemplified by an annulus bounded by a pair of concentric circles in
$\RR^2$. The inner boundary is the circle of smaller radius and has
total curvature $-2\pi$, while the outer boundary has total curvature
$2\pi$.

A regular annulus is {\em primitive} with respect to $q$ if it
contains no singularities of $q$ in its interior. By
Equation~(\ref{eqn:gaussbonnet}), its boundaries satisfy
$\kappa_A(\dd_0) = - \kappa_A(\dd_1)$. A regular annulus is {\em flat}
if $\kappa_A(\dd_0)=\kappa_A(\dd_1)=0.$ By~(\ref{eqn:gaussbonnet}), a
flat annulus is necessarily primitive, and is foliated by Euclidean
geodesics homotopic to the boundaries.  Thus a flat annulus is
isometric to a cylinder obtained as the quotient of a Euclidean
rectangle in $\RR^2$.  Note that a primitive annulus must either be
flat or expanding.

One reason for introducing flat and expanding annuli is that their
moduli are easy to estimate.  The following result can be deduced from
Theorem~4.5 of \cite{minskyharmonic} and is proved in \cite{rafi1}:
\begin{theorem}
\label{thm:modcomparison} Let $A \subset S$ be an annulus
that is primitive with respect to $q$ and with inner and outer
boundaries $\dd_0$ and $\dd_1$, respectively.  Let $d$ be the
$q$-distance between $\dd_0$ and $\dd_1$. Then either
\begin{enumerate}
\item[(i)] $A$ is flat and $\mod A = d/l_q(\dd_0)$ or
\item[(ii)] $A$ is expanding and $\mod A \asymp \log [d/l_q(\dd_0)]$.
\end{enumerate}
\end{theorem}

\subsection{Modulus of annulus and length of short curve}
\label{sec:modofann}
The link between the hyperbolic and quadratic differential metrics on
a surface is made using annuli of large modulus.  Let $\s$ be the
hyperbolic metric that uniformizes $\Sigma$.  If $\a$ is short in
$\s$, Maskit \cite{maskit} showed that the extremal length
$\Ext_\Sigma(\a)$ and hyperbolic length $l_\s(\a)$ are comparable, up
to multiplicative constants. Moreover, there is an embedded collar
$C(\a)$ around $\a$ in $(S,\s)$ whose modulus is comparable to
$1/l_\s(\a)$ (see \cite{minskyproduct} for an explicit calculation),
and therefore also to $1/\Ext_\Sigma(\a)$.  By the (geometric)
definition of extremal length, this implies that the maximal annulus
around $\a$ in $\Sigma$ has modulus comparable to $1/l_\s(\a)$. The
following theorem of Minsky allows us to replace any annulus of
sufficiently large modulus with one that is primitive:
\begin{theorem} (\cite{minskyharmonic} Theorem 4.6)
\label{thm:prim} Let $A \subset \Sigma$ be any homotopically non-trivial 
annulus whose modulus is sufficiently large and let $q \in {\mathcal
  Q}(\Sigma)$. Then $A$ contains an annulus $B$ that is primitive with
respect to $q$ and such that $\mod A \asymp \mod B$.
\end{theorem}
\noindent (The statement of Theorem 4.6 in
\cite{minskyharmonic}  should read $\mod A \geq m_0$ not $\mod A \leq
m_0$.)
 Thus, we have
\begin{theorem}
\label{basic}
If $\a$ is a simple closed curve which is sufficiently short in
$(S,\s)$, then for any $q \in {\cal Q}(\Sigma)$, there is an annulus
$A$ that is primitive with respect to $q$ with core homotopic to $\a$
such that
$$\frac{1}{l_\s(\a)} \asymp \mod (A).$$  
\end{theorem}
\noindent We may assume that $\ep_0$ was chosen so that if 
$l_\s(\a) < \ep_0$, then $l_\s(\a)$ is small enough that this theorem is
valid.

We can now apply Theorem~\ref{thm:modcomparison} in the following way.
It follows from Equation~(\ref{eqn:gaussbonnet}) that every simple
closed curve $\gamma$ on $(S,q)$ either has a unique $q$-geodesic
representative, or is contained in a family of closed Euclidean
geodesics which foliate a flat annulus~\cite{strebel}. Denote by
$F(\gamma)$ the maximal flat annulus, which necessarily contains all
$q$-geodesic representatives of $\gamma$.  If the geodesic
representative of $\gamma$ is unique, then $F(\gamma)$ is taken to be
the degenerate annulus containing this geodesic alone.  Denote the
(possibly coincident) boundary curves of $F(\gamma)$ by $\dd_0, \dd_1$
and consider the $q$-equidistant curves from $\dd_i$ outside
$F(\gamma)$.  Let $\hat {\dd_i}$ denote the first such curve which is
not embedded.  If $\hat \dd_i \neq \dd_i$, then the pair $\dd_i$,
$\hat {\dd_i}$ bounds a region $ {E_i}(\gamma)$ whose interior is an
annulus with core homotopic to $\gamma$, and which is by its
construction regular and expanding.  Combining the preceding two
theorems with Theorem~\ref{thm:modcomparison} we have:
\begin{corollary}
\label{cor:lengthest}  If $\alpha$ is an extremely short curve 
on $(S,\s)$, then
$$
\frac{1}{l_{\s}(\alpha)} \asymp 
\max \left\{\mod F(\alpha), \mod {E_0}(\alpha),\mod {E_1}(\alpha) \right\}.
$$
\end{corollary}
\begin{proof}
  Since $\a$ is extremely short, we have $1/\Ext_\Sigma(\a) \asymp
  1/l_\s(\a)$ and hence it follows from the (geometric) definition of
  extremal length that
$$
\frac{1}{l_\s(\a)} \asymp \frac{1}{\Ext_\Sigma(\a)} \geq 
\max \left\{ \mod F(\alpha), \mod {E_0}(\alpha),\mod {E_1}(\alpha) \right\}.
$$

By Theorem~\ref{basic}, there is a primitive annulus $A$ whose core is
homotopic to $\a$ such that $1/l_\s(\a) \asymp \mod (A).$ We will show
that
$$\mod (A) \prec \max \{\mod F(\alpha), \mod {E_0}(\alpha), \mod
{E_1}(\alpha) \}.$$
If $A$ is flat, then $A$ must be contained in the
maximal flat annulus $F(\a)$. In this case, $\mod (A) \leq \mod
F(\a)$. If $A$ is expanding, then although it may not be contained in
either $E_0(\a)$ or $E_1(\a)$ it must be disjoint from the interior of
$F(\a)$. Without loss of generality, let us assume that $A$ lies on
the same side of $F(\a)$ as $E_0(\a)$.  Let $\dd_0$ and $\hat \dd_0$
be respectively, the inner and outer boundaries of $E_0(\a)$. Let
$C_0$ and $C_1$ be respectively, the inner and outer boundaries of
$A$.  Since $l_{q}(\dd_0)$ is equal to the $q$-length of the geodesic
representative of $\a$, the $q$-length of the inner boundary of $A$
satisfies $l_q(C_0) \geq l_q(\dd_0)$. Let $\omega$ be a $q$-shortest
arc in $E_0(\a)$ from $\dd_0$ to itself; its length is $2d_q(\dd_0,
\hat \dd_0)$. The intersection $\omega \cap A$ is a union of two arcs,
each of which goes from one boundary component of $A$ to another.
Since $l_q(\omega \cap A) \leq l_q(\omega)$, it follows that
$d_q(C_0,C_1) \leq l_q(\omega \cap A)/2 \leq d_q(\dd_0,\hat \dd_0)$.
Thus, it follows from Theorem~\ref{thm:modcomparison} that
\begin{equation*}
\mod (A) \asymp \log \frac{d_q(C_0,C_1)}{l_q(C_0)} \leq  
\log \frac{d_q(\dd_0,\hat \dd_0)}{l_q(\dd_0)} \asymp \mod (E_0(\a)). \qedhere
\end{equation*}
\end{proof}
The idea of our basic length estimates for extremely short curves
$\alpha$ in Theorem~\ref{prop:shortingeo}, is to combine this
corollary with the estimates for the moduli of $F(\alpha)$ and $
{E_i}(\alpha)$ in Theorem~\ref{thm:modcomparison}.

\subsection{Thick-thin decomposition and the $q$-metric}
\label{sec:thickthin}
The thick-thin decomposition for quadratic differentials developed
in~\cite{rafi3} describes the relationship between the $q$-metric on
the surface $\Sigma$ and the uniformizing hyperbolic metric $\s$ in
the thick components of the thick-thin decomposition of $\s$. It
states that on the hyperbolic thick parts of $(S,\s)$ the two metrics
are comparable, up to a factor which depends on the moduli of the
expanding annuli around the short curves in the boundary of the thick
component.  This factor will be crucial in our estimates below.

To make a precise statement, for a subsurface $Y$ of $S$, define the
$q$-{\em geodesic representative of} $Y$ to be the unique subsurface
$\hat Y$ of $(S,q)$ with $q$-geodesic boundary in the homotopy class
of $Y$ that is disjoint from the interior of $F(\gamma_i)$ for all
components $\gamma_i $ of $ \dd Y$.  Notice that $\hat Y$ is
$q$-geodesically convex, so that if a closed curve $\zeta$ is
contained in $Y$, it has a $q$-geodesic representative contained in
$\hat Y$. (It is possible for $\hat Y$ to be degenerate. See
\cite{rafi3} for an example where the area of $\hat Y$ is zero.)

If $Y$ is not a pair of pants, define $\lambda_Y$ to be the length of
the $q$-shortest non-peripheral simple closed curve contained in $\hat
Y$.  If $Y$ is a pair of pants, define $\lambda_Y$ to be $\max \{
l_q(\gamma_1),l_q(\gamma_2),l_q(\gamma_3)\}$ where $\gamma_1,
\gamma_2,\gamma_3$ are the three boundary curves of $\hat Y$. The
thick-thin decomposition for quadratic differentials is the following:
\begin{theorem} [Rafi \cite{rafi3}] 
\label{thm:scalefactor} Let $\s$ be the hyperbolic metric that
uniformizes $\Sigma$ and let $Y$ be a thick component of the
hyperbolic thick-thin decomposition of $(S,\s)$.  Then
\begin{enumerate}
\item[(i)] $diam_q \hat Y \mul \lambda_Y$,
\item[(ii)] For any non-peripheral simple closed curve $\zeta$ in $  Y$, 
we have
$$l_{q}(\zeta) \mul  {\lambda_Y} l_{\s}(\zeta).$$
\end{enumerate}
\end{theorem}

\subsection{Twist in the $q$-metric}
\label{sec:qtwist}
 In order to compare two surfaces, we need to estimate not only the
lengths  but also the twist parameters of short curves. To do this we 
use a signed version
of Rafi's definition (\cite{rafi2} Section 4)  of the twist 
of a simple curve $\zeta$  about another  curve $\a$ in a
quadratic differential metric $q$ on $S$.

Let $A \subset S$ be a regular annulus with core curve $\a$, let
$\tilde S$ be the annular cover of $S$ corresponding to $\a$, and let
$\tilde \a$ be a $\tilde q$-geodesic representative of the core of
$\tilde S$.  Suppose that $\zeta$ is a simple $q$-geodesic (i.e.
geodesic with respect to the $q$-metric) that intersects $\a$, and let
$\tilde \zeta$ be a lift which intersects $\tilde \a$.  Let $\tilde
\beta$ be a bi-infinite $\tilde q$-geodesic arc in $\tilde S$ that is
orthogonal to $\tilde \a$.  We would like to define the twist
$tw_q(\zeta,\a)$ to be the sum $a_{\tilde S}(\tilde \zeta, \tilde
\beta)$ of the signed intersection numbers over all intersections
between $\tilde \zeta$ and $\tilde \beta$.  The following lemma shows
that $a_{\tilde S}(\tilde \zeta, \tilde \beta)$ is, up to a bounded
additive error, independent of the choices of $\tilde \beta$ and
$\tilde \zeta$ .

\begin{lemma}
\label{lem:flatintersection} Let $\alpha$, $\tilde \alpha$, and $\tilde S$  be
as above and suppose that $\zeta $ is a simple $q$-geodesic transverse
to $\alpha$.  Suppose that $\tilde \zeta, \tilde \zeta'$ are different
lifts of $\zeta$ that intersect $\tilde \alpha$ and that $ \tilde
\beta, \tilde \beta'$ are different bi-infinite $\tilde q$-geodesic
arcs orthogonal to $\tilde \alpha$. Then
$$
a_{\tilde S}(\tilde \zeta, \tilde \beta)=
a_{\tilde S}(\tilde \zeta', \tilde \beta') \pm O(1).
$$
\end{lemma} 
\begin{proof} 
  Let $\tilde F(\a)$ be the lift to $\tilde S$ of the maximal flat
  annulus $ F(\a)$ around $\tilde \a$ and let $a_{ F}(\tilde \zeta,
  \tilde \beta)$ denote the sum of signed intersection numbers over
  intersection points within $\tilde F(\a)$. A simple Euclidean
  argument shows that for any two disjoint arcs $\tilde \zeta, \tilde
  \zeta'$, we have
$$
a_{F}(\tilde \zeta, \tilde \beta)=a_{ F}(\tilde \zeta', \tilde \beta') \pm O(1).
$$
We claim that outside $\tilde F(\a)$, any two $\tilde q$-geodesics
can intersect at most twice. Outside $\tilde F(\a)$, $\tilde S$ is
made up of two regular expanding annuli $E_1,E_2$, one attached to
each boundary of $\tilde F(\a)$.  These annuli extend out to infinity
in $\tilde S$ (which can be compactified using the hyperbolic metric
on $S$, see~\cite{masurminskyII} Section 2.4). The key point is that
in any expanding annulus $E$, two geodesic arcs can intersect at most
once. For if they intersected twice, we would get a piecewise geodesic
loop $\gamma$ homotopic to the inner boundary, made up of two geodesic
arcs that go from one intersection point to the other. Along each arc,
the geodesic curvature vanishes. The Gauss-Bonnet theorem in
Equation~(\ref{eqn:gaussbonnet}) applied to the annulus bounded by
$\gamma$ and the inner boundary shows this is impossible.  (Notice
that if $E$ is not primitive, then the singularities of $q$ in $E$
only improve the desired inequality in
Equation~(\ref{eqn:gaussbonnet}).)
\end{proof}

We define $tw_q(\zeta,\a)$ to be the minimum of the numbers $a_{\tilde
  S}(\tilde \zeta, \tilde \beta) $ over all choices of $\tilde \zeta,
\tilde \beta$. Notice that the argument requires only that the lifts
of $\zeta$ to $\tilde S$ be disjoint, so that we can similarly define
$tw_q(\nu,\a)$ for a geodesic lamination $\nu$ where as usual,
$tw_q(\nu,\a)$ depends only the underlying support $|\nu|$ of $\nu$.

The following key result allows us directly to compare the twists in
the hyperbolic and quadratic differential metrics.
\begin{proposition}[\cite{rafi2} Theorem 4.3] 
 \label{prop:twistcomparison} Suppose that
 $\s$ is a hyperbolic metric uniformizing a surface $\Sigma \in \ts$
 and that $q \in \Q(\Sigma )$, and let $\nu$ be a geodesic lamination
 intersecting $\a$. Then
 $$
 |tw_{\s}(\nu,\a) - tw_q(\nu,\a)| = O( 1/l_{\s}(\a)).$$
\end{proposition}
The statement in~\cite{rafi2} has $ i_{\a}(\nu,\delta_{\a})$ in place
of $ tw_{\s}(\nu,\a)$, however by Lemma~\ref{lemma:twistintersect},
this distinction is unimportant.  The result in~\cite{rafi2} is stated
for closed curves but extends immediately to the case of geodesic
laminations as explained above.

\subsection{Length and twist along $\G(\np,\nm)$.}

As explained at the end of Section~\ref{sec:modofann}, we can use
Theorem~\ref{thm:modcomparison} and Corollary~\ref{cor:lengthest} to
estimate the length of an extremely short curve $\a$ in $\G_t$.  We
call a flat or expanding annulus which achieves the maximum modulus in
Corollary \ref{cor:lengthest} a \emph{dominant annulus} for $\a$.
There may be more than one dominant annulus, but this will not affect
our reasoning and we will refer to `the dominant annulus'. The
estimates depend on whether the dominant annulus is flat or expanding,
corresponding to the two terms $D_t(\a)$ and $K_t(\a)$ in the main
result Theorem~\ref{prop:shortingeo} of this section.

Suppose first that the maximal flat annulus $ F_t(\a)$ is dominant.
Provided that $\a$ is neither vertical nor horizontal (see
Section~\ref{sec:balancetime}), the following proposition expresses
$\mod F_t(\a)$ in terms of the relative twisting $d_\a(\np,\nm)$ of
$\np,\nm$ around $\a$ defined in Section~\ref{sec:reltwist}.  The case
in which $\a$ is either horizontal or vertical, so that either $\np$
or $\nm$ has empty intersection with $\a$, is easier and is dealt with
in Section~\ref{sec:vertonG}.

\begin{proposition}
  Let $\a$ be a curve in $(S,q)$ that is neither vertical nor
  horizontal and suppose that the maximal flat annulus $ F_t(\a)$ is
  dominant.  Then \label{prop:flat}
  $$
  \mod F_t(\a) \asymp e^{-2|t-t_\a|} d_\a(\np,\nm).$$
\end{proposition}

 \begin{proof}
   Since a flat annulus is Euclidean, its geometry is very simple.
   Let $\eta$ be a $q_t$-geodesic arc in $F_t(\a)$ joining the two
   boundaries of $F_t(\a)$ that is orthogonal to the geodesic
   representatives of $\a$.  For a simple geodesic $\zeta$ transverse
   to $\a$, define $tw_{F_t}(\zeta,\a)$ to be the signed intersection
   number of $\zeta$ with $\eta$ in $F_t(\a)$.  It is independent of
   the choice of $\eta$ up to a bounded error of $1$.  Assuming that
   $\a$ is neither vertical nor horizontal, then at the balance time
   $t_\a$ (see Section~\ref{sec:balancetime}) the horizontal and
   vertical foliations both make an angle of $\pi/4$ with the
   $q_{t_\a}$-geodesic representatives of $\alpha$. In this case, a
   leaf of $\nu_{t_\a}^+$ or $\nu_{t_\a}^-$ intersects $\eta$
   approximately (up to an error of $1$) $l_{q_{t_\a}}(\eta) /
   l_{q_{t_\a}}(\a)$ times, so the modulus of $F_{t_\a}(\a)$ is
   approximated by $tw_{F_{t_\a}}(\np,\a)= tw_{F_{t_\a}}(\nm,\a)$.
   More generally, the horizontal leaves make an angle $\psi_t$ with
   $\alpha$, where $|\tan \psi_t| = e^{2(t-t_\a)}$. From this it is a
   straightforward exercise in Euclidean geometry, see Section 4.1
   of~\cite{rafi2}, to prove:
\begin{equation}
\label{eqn:geomflatann}
|tw_{F_t}(\nu^{\pm},\a ) - e^{{\mp}2(t-t_\a)} \mod F_t(\a)|  \leq 1.
\end{equation}

We will show that 
\begin{equation}
\label{eqn:twistinann} 
|tw_{F_t}(\np,\a) - tw_{F_t}(\nm,\a)| = d_\a(\np,\nm) \pm O(1),
\end{equation}
from which the proposition follows.  From the proof of
Lemma~\ref{lem:flatintersection} we have $tw_{q_t}(\zeta,\a) =
tw_{F_t}(\zeta,\a) \pm O(1)$.  Now it was observed in
Section~\ref{sec:fncoords} (see also the subsequent discussion in
Section~\ref{sec:reltwist}) that although $tw_\s(\nu,\a)$ depends on
the metric $\s$ in which it is measured, the difference in twist of
$\np$ and $\nm$ equals (up to a bounded error) the number of times a
leaf of $\np$ intersects a leaf of $\nm$ in the annular cover of $S$
corresponding to $\a$. Since this also holds for a quadratic
differential metric, we get:
\begin{align*}
|tw_{q_t}(\np,\a) - tw_{q_t}(\nm,\a)| &= |tw_\s(\np,\a) - tw_\s(\nm,\a)| 
\pm O(1) \\
&= \da \pm O(1).
\end{align*}
Equation~(\ref{eqn:twistinann})   follows.
\end{proof}
 
\medskip

Suppose now that one or other of the expanding annuli around $\a$ is
dominant.  The estimate of modulus in this case is given by:
\begin{proposition}
\label{expanding} Let $q \in {\cal Q}(\Sigma)$. 
Suppose that $\a$ is extremely short in $\s$ and let $Y$ be a thick
component of the hyperbolic thick-thin decomposition of $(S,\s)$, one
of whose boundary components is $\a$. Let $\hat \a$ be the
$q$-geodesic representative of $\a$ on the boundary of $\hat Y $ and
let $E(\a)$ be a maximal expanding annulus on the same side of $\hat
\a$ as $\hat Y$. If $E(\a)$ is dominant, then
$$\mod E(\a) \asymp \log \frac{\lambda_Y}{l_q(\a)}.$$
\end{proposition}
\begin{proof} Let $d_q$ denote the $q$-metric. By 
  Theorems~\ref{thm:prim} and~\ref{thm:modcomparison}(ii), it is
  sufficient to show that $d_q(\dd_0,\dd_1) \asymp \lambda_Y$.
  
  Note that although $E(\a)$ is not necessarily contained in $\hat Y$,
  the outer boundary $\dd_1$ must intersect $\hat Y$. Hence,
  $d_q(\dd_0,\dd_1) \leq \mbox{diam}_q (\hat Y)$ and so by
  Theorem~\ref{thm:scalefactor}(i), we have $d_q(\dd_0,\dd_1) \prec
  \lambda_Y.$
  
  Now we prove the inequality in the other direction. Observe that
  since $E(\a)$ is maximal, the outer boundary $\dd_1$ intersects
  itself and so there is a non-trivial arc $\omega$ with endpoints on
  $\hat \a$ whose length is $2d_q(\dd_0, \dd_1)$. First suppose that
  $\omega$ is contained in $\hat Y$. A regular neighborhood of $\hat
  \a \cup \omega$ is a pair of pants whose boundary curves are
  homotopic to $\a$ and two additional curves $\zeta_1,\zeta_2$. Note
  that for $i=1,2$,
  $$l_q(\omega) + l_q(\a) \geq l_q(\zeta_i).$$
  Thus, if either $
  \zeta_1$ or $\zeta_2$, say $\zeta_1$ is non-peripheral in $Y$, then
  $$
  \frac{d_q(\dd_0,\dd_1)}{l_q(\a)} \succ
  \frac{l_q(\zeta_1)}{l_q(\a)} \geq \frac{\lambda_Y}{l_q(\a)}.$$
  If
  both $\zeta_1,\zeta_2$ are peripheral, then $Y$ is a pair of pants,
  and we have
$$\frac{d_q(\dd_0,\dd_1)}{l_q(\a)} \succ \frac{l_q(\zeta_1) +
  l_q(\zeta_2)}{l_q(\a)} \asymp \frac{\lambda_Y}{l_q(\a)}.$$

If $\omega$ exits $\hat Y$, we replace it with a new arc $\omega'$ as
follows. Let $p$ be the first exit point and let $\gamma$ be the
boundary component of $\hat Y$ that contains $p$.  Let $\omega'$ be
the arc that first goes along $\omega$ to $p$, then makes one turn
around $\gamma$ from $p$ to itself, then comes back to $\hat \a$ along
the first path. Because $\gamma$ is in the boundary of $\hat Y$, its
hyperbolic length is extremely short and thus by Theorem~\ref{basic},
the original arc $\omega$ must pass through an annulus of large
modulus with core curve $\gamma$. Therefore, $l_q(\gamma) \prec
l_q(\omega)$ and so we have
$$l_q(\omega') \leq 2l_q(\omega) + l_q(\gamma) \prec l_q(\omega).$$
Now we can run the same argument as above with $\omega'$ in place of
$\omega$ to deduce the desired inequality.
\end{proof}

We are now able to write down the desired length estimate. Suppose
that the curve $\a$ is extremely short in some surface $\mathcal G_t$
along the \Teich geodesic $\G(\np,\nm)$.  Let $Y_1,Y_2$ be the thick
components of the thick-thin decomposition of $(S,\mathcal G_t)$ that
are adjacent to $\a$ (where $Y_1$ may equal $Y_2$).  Define
\begin{align}
\label{ktalpha}
K_t(\a) &= \max \left\{ \frac{\lambda_{Y_1}}{l_{q_t}(\a)},
\frac{\lambda_{Y_2}}{l_{q_t}(\a)} \right\}  
\intertext{and}
\label{dtalpha}
D_t(\a) &=e^{-2|t-t_\a|} d_\a(\np,\nm).
\end{align}
Then, combining Corollary~\ref{cor:lengthest} and
Propositions~\ref{prop:flat} and~\ref{expanding}, we obtain our first
main result which is essentially Theorem~\ref{shortong} of the
Introduction:
\begin{theorem}
\label{prop:shortingeo}
Let $\a$ be a curve on $S$ that is neither vertical nor horizontal.
If $\a$ is extremely short in $\mathcal G_t$, then
$$\frac{1}{l_{\mathcal G_t}(\a)} \asymp \max\{ D_t(\a), \log K_{t}(\a) \}.$$
\end{theorem}

It was shown in \cite{rafi1} that $K_t(\a)$ can be estimated
combinatorially as the subsurface intersection of $\nu^+,\nu^-$ in the
corresponding component of the thick part of $\mathcal G_t$ adjacent
$\a$.  Since this is not
necessary for our development, we shall not go into this here.\\

We need to estimate not only the lengths but also the twist parameters
about short curves.  Combining Proposition~\ref{prop:twistcomparison}
with Equation~(\ref{eqn:geomflatann}),
Lemma~\ref{lem:flatintersection}, Proposition~\ref{prop:flat}, and
Theorem~\ref{prop:shortingeo}, we get:
\begin{theorem}[\cite{rafi2} Theorem 1.3]
\label{prop:quadtwistest}
With the same hypotheses as Theorem~\ref{prop:shortingeo}, 
\begin{align*}
Tw_{\mathcal G_t}(\np,\a) &=
O \left(\frac{1}{l_{\mathcal G_t}(\a)} \right) \ \text{if} 
\quad t \geq t_\a,  \\
Tw_{\mathcal G_t}(\nm,\a) &=
O \left(\frac{1}{l_{\mathcal G_t}(\a)} \right) \ \text{if} \quad t \leq t_\a.
\end{align*}
\end{theorem}

\subsection{Estimates for horizontal and vertical short curves on $\G$}
\label{sec:vertonG}
We also need the analogue of Theorems~\ref{prop:shortingeo} and
\ref{prop:quadtwistest} for extremely short curves $\a$ which are
either horizontal or vertical.

For definiteness, assume $\a$ is vertical so that $i(\a,\nu^-) = 0$.
The definition of balance time no longer makes sense.  Instead, we
work relative to time $t=0$.  Let $d$ be the height (i.e., distance
between the two boundaries) of $F_{0}(\a)$.  At an arbitrary time $t$,
$$l_{q_{t}}(\a) = i(\a,\nu^+_t)= e^{t}l_{q_{0}}(\a)$$ 
while the height of $F_t(\a)$ is $e^{-t} {d}$. Hence
\begin{equation}
\label{eqn:flat1}
\mod F_{t}(\a) = e^{-2t}\mod F_{0}(\a) .
\end{equation}
This is the analogue of Proposition~\ref{prop:flat}.
  
If $\a$ is vertical, the discussion in Proposition~\ref{expanding}
about expanding annuli is unchanged. Thus we obtain:
\begin{theorem}
\label{prop:shortingeo1}
Let $\a$ be a vertical curve on $S$.  If $\a$ is extremely short in
$\mathcal G_t$, then
$$
\frac{1}{l_{\mathcal G_t}(\a)} \asymp 
\max \left\{ e^{-2t}\mod F_{0}(\a), \log K_{t}(\a) \right\}.
$$ 
If $\a$  is  horizontal, the estimate is the same except
that the first term is replaced by  $e^{2t}\mod F_{0}(\a)$.
\end{theorem}

We also want the analogue of Theorem~\ref{prop:quadtwistest}.  
If $\a$ is vertical, then $Tw_{\mathcal G_t}(\nm,\a) $ is undefined. 
However we  have:
\begin{theorem} 
\label{prop:quadtwistest1}
If $\a$ is extremely short in $\mathcal G_t$, and if $\a$ is
vertical then 
\begin{align*} 
Tw_{\mathcal G_t}(\np,\a) & = O \left(\frac{1}{l_{\mathcal G_t}(\a)} \right), \\
\intertext{while if $\a$ is horizontal then}
Tw_{\mathcal G_t}(\nm,\a) & = O \left(\frac{1}{l_{\mathcal G_t}(\a)} \right). 
\end{align*}
\end{theorem}
\begin{proof}
  If $\a$ is vertical, the $q$-twist $tw_q(\nu^+, \a)$ in $F_t(\a)$
  vanishes, while if $\a$ is horizontal, then $tw_q(\nu^-, \a)=0$.
  The result follows from Lemma~\ref{lem:flatintersection} and
  Proposition~\ref{prop:twistcomparison}.
\end{proof}
\section{Short curves along  lines of minima}
\label{sec:shortlm}
In this section we prove Theorem B, stated more precisely as
Theorem~\ref{thm:shortonL}, which gives our combinatorial estimate for
the length of a curve which becomes extremely short at some point
along the line of minima $\L(\nu^+,\nu^-)$. We also estimate the twist
of $\nu^\pm$ around $\a$ in Theorem~\ref{thm:twistonL}.  This will
form the basis for our comparison of the metrics $\L_t$ and $\mathcal
G_t$.  It turns out that, in a close parallel to the case of the
Teichm\"uller geodesic, there are two reasons why a curve can be
extremely short: either the relative twisting of $\nu^+,\nu^-$ about
$\a$ is large, or one or other of the pants curves in a pair of pants
adjacent to $\a$ in a short marking of $\L_t$ has large intersection
with either $\nu^+$ or $\nu^-$.

More precisely, suppose that $\a$ is an extremely short curve in
$\L_t$ and let ${\cal P}_{\L_t}$ be a short pants system in $\L_t$,
which necessarily contains $\a$. Define
\begin{equation}
\label{eqn:defofH}
H_t(\a)=  
\sup_{\beta \in {\cal B}} \frac{l_{q_t}(\beta)}{l_{q_t}(\a)}, 
\end{equation}
where $\cal B$ is the set of pants curves in ${\cal P}_{\L_t}$ which
are boundaries of pants adjacent to $\a$ and $q_t$ is the quadratic
differential metric of area $1$ (on the corresponding surface $\G_t$)
whose horizontal and vertical foliations are respectively, $\nu_t^+$
and $\nu_t^-$.

Let $D_t(\a)$ be as in Equation~(\ref{dtalpha}).  Our main estimates are:
\begin{theorem}
\label{thm:shortonL}
Let $\a$ be a curve on $S$ which is neither vertical nor horizontal.
If $\a$ is extremely short in $\L_t$, then
$$
\frac{1}{l_{\L_t}(\a)} {\asymp} \max
\left\{ D_t(\a), \sqrt{H_t(\a)} \right\}. 
$$
\end{theorem}

\begin{theorem}
\label{thm:twistonL} 
With the same hypotheses as Theorem~\ref{thm:shortonL}, the twist
satisfies:
\begin{align*}
Tw_{\L_t}(\nu^+,\a) 
 &= O\left( \frac{1}{ l_{\L_t}(\a)} \right) \ \text{if} \quad t \geq t_\a , \\
Tw_{\L_t}(\nu^-,\a) 
&= O\left( \frac{1}{ l_{\L_t}(\a)} \right) \ \text{if} \quad t \leq t_\a.
\end{align*} 
\end{theorem}

To prove Theorems~\ref{thm:shortonL},~\ref{thm:twistonL}, we note that
since the surface $\L_t$ is on the line of minima, we have at the
point $\L_t$,
\begin{equation}
\label{eqn:1form}
\mbox{d} l(\nu^+_t) + \mbox{d} l(\nu^-_t) = 0.
\end{equation}
The pants curves in ${\mathcal P}_{\L_t}$ (together with seams) define
a set of coordinates $(l_\s(\a),\tau_\a(\s))$ on $\ts$ as explained in
Section~\ref{sec:fncoords}, which in turn define infinitesimal twist
$\frac{\dd}{\dd {\tau}_\a}$ and length $\frac{\dd}{\dd l(\a)}$
deformations for $\a \in {\mathcal P}_{\L_t}$.
Theorems~\ref{thm:shortonL} and~\ref{thm:twistonL} will follow from
the relations we get by applying Equation~(\ref{eqn:1form}) to
$\frac{\dd}{\dd {\tau}_\a}$ and $\frac{\dd}{\dd l(\a)}$.  For
$\frac{\dd}{\dd {\tau}_\a}$, we use the well-known formula of
Kerckhoff~\cite{kercktwist} and Wolpert~\cite{wolpert}, while for
$\frac{\dd}{\dd l(\a)}$ we use the analogous formula for the length
deformation derived in~\cite{cmswolpert}.

\subsection{Differentiation with respect to twist}
\label{sec:kerckhoff}
Suppose as above that $\a$ is an extremely short curve in $\L_t$.  If
we apply Equation~(\ref{eqn:1form}) to $\frac{\dd}{\dd \tau_\a}$, the
derivative formula in \cite{kercktwist} and \cite{wolpert} gives
$$
0=\frac{\dd l(\nu^+_t)}{\dd \tau_\a} + 
\frac{\dd l(\nu^-_t)}{\dd \tau_\a}=
\int_\a \cos \theta^+ d\nu_t^+ + \int_\a \cos \theta^- d\nu_t^-,
$$
where $\theta^\pm$ is the function measuring the angle from each
arc of $|\nu^\pm_t|$ to $\a$. Assume that $\a$ is neither vertical nor
horizontal, so that neither $i(\nu^+,\a)$ nor $i(\nu^-,\a)$ is zero.
Then we may define the average angle $\Theta_t^{\pm}$ by
$$
\cos \Theta_t^{\pm}
= \frac{1}{i(\nu^\pm_t,\a)} \int_\a \cos \theta^\pm d\nu_t^\pm.
$$
Setting $T = t-t_{\a}$, the preceding two equations give
\begin{equation}
\label{eqn:kerckhoff}
e^{T}\cos \Theta^+_t + e^{-T}\cos \Theta^-_t = 0.
\end{equation}

If a particular leaf $L$ of a lamination $|\nu|$ cuts $\a$ at an angle
$\theta$ at a point $p$, then from the definition of the twist (see
Section~\ref{sec:twistshyperbolic}) and simple hyperbolic geometry we
have
$$
\cos \theta  = \tanh \frac{ tw_{\L_t} (L ,\a,p) \, l_{\L_t}(\a)}2.
$$
Since the twists $ tw_{\L_t} (L ,\a,p)$ for different leaves
$L$ differ by at most $1$, if $\alpha$ is sufficiently short we
obtain the estimate 
\begin{equation}
\label{eqn:twistangle}
|\cos \theta - \cos \Theta_t^{\pm} | = O(l_{\L_t}(\a)) 
\end{equation}
from which we deduce that either $\cos \Theta_t^{\pm}$ and
$tw_{\L_t}(\nu^\pm,\a)$ have the same sign, or that $|\cos
\Theta_t^{\pm} | = O(l_{\L_t}(\a))$ so that $|tw_{\L_t}(\nu^\pm,\a)| =
O(1)$.

Note also that Equation~(\ref{eqn:kerckhoff}) implies that $\np,\nm$
twist around $\a$ in opposite directions and that the lamination whose
weight on $\a$ is smaller does more of the twisting.
\subsection{Differentiation with respect to length}
\label{sec:differentiation}
For the length deformation, we shall apply the extension of the
Wolpert formula derived in \cite{cmswolpert}, which gives a general
expression for $\mbox{d} l(\zeta)$, for $\zeta \in \S$, with reference
to a pants curves system $\mathcal P$.  Let $\tilde \a_1, \ldots,
\tilde \a_n$ be the lifts of the pants curves in $\cal P$ successively
met by $\zeta$, where the segment of the lift $\tilde \zeta$ of
$\zeta$ between $\tilde \a_1$ and $\tilde \a_n$ projects to one
complete period of $\zeta$. Let $d_j$ be the length of the common
perpendicular $\pi_j$ between $\tilde \a_j, \tilde \a_{j+1}$ and let
$S_j$ be the signed distance between $\pi_{j-1}$ and $\pi_j$ along
$\tilde \a_j$, where the sign is positive if the direction from
$\pi_{j-1}$ to $\pi_j$ coincides with the orientation of $\a_j$. (Note
that if $\tilde \a_j$ and $\tilde \alpha_{j+1}$ project to the same
curve $\alpha$ and if $\alpha$ is adjacent to two distinct pairs of
pants, then $\pi_j$ projects to an arc perpendicular to $\alpha$ that
is not a seam.)  Then Equation~(3) of \cite{cmswolpert} states that
\begin{equation}
\label{eqn:cmswolpert}
\mbox{d}l(\zeta)=\sum_{j=1}^n \cosh u_j \,\mbox{d}d_j +
  \sum_{j=1}^n \cos \theta_j \,\mbox{d} S_j,
\end{equation}
where $\theta_j$ is the angle from $\tilde \zeta$ to $\tilde \a_j$
measured counter-clockwise and $u_j$ is the complex distance from
$\tilde \zeta$ to the complete bi-infinite geodesic which contains
$\pi_j$. Replacing sums by integrals, we see that this formula,
derived in~\cite{cmswolpert} for closed curves, pertains equally to a
measured lamination.
  
In our case, we take $\cal P$ to be $\cal P_{\L_t}$ and apply this
formula to $\frac{\dd}{\dd l(\a)}$ for $\a \in \cal P_{\L_t}$.  The
non-zero contributions will be from terms $\mbox{d}S_j$ corresponding
to lifts of $\a$, and from two types of terms $\mbox{d}d_j$: those
corresponding to perpendiculars with endpoints on lifts of $\a$, and
those corresponding to perpendiculars which do not intersect any lift
of $\a$, but whose projections are contained in a common pair of pants
with $\a$.
  
We first estimate the contribution from the terms $\mbox{d}d_j$.
Suppose as above that ${\cal P_{\L_t}}$ is a short pants decomposition
for $\cal L_t$. Let $P$ be a pair of pants in $S\setminus{\cal
  P_{\L_t}}$ that has $\a$ as a boundary component.  The geometry of
$P$ is completely determined by the lengths of the three boundary
curves $\a,\beta,\gamma$.  A common perpendicular joining two (not
necessarily distinct) boundary components of $P$ may or may not have
one of its endpoints on a boundary curve which projects to $\alpha$ on
$S$. We say that the common perpendiculars of the first kind are
\emph{adjacent} to $\alpha$, while those of the second type are not.
The terms $\mbox{d}d_j$ are estimated by the following lemma which is
proved in the Appendix:
\begin{lemma}
\label{lem:hexagon*}
Suppose that $\a$ is extremely short in $\L_t$ and let
$P$ be a pair of pants in $S\setminus{\cal
  P_{\L_t}}$ that has $\a$ as a boundary component.
Let $v$ denote the length of a
common perpendicular adjacent to $\a$, and let $w$ denote the length
of a common perpendicular not adjacent to $\a$. Then
$$
\frac {\dd v}{\dd  l(\a)} \mul -\frac{1}{l(\a)}
\quad\text{and}\quad
\frac{\dd w}{\dd l(\a)}  \mul l(\a),
$$
where the partial derivatives are taken with respect to the
coordinates $(l(\a),\tau_\a)_{\a \in {\cal P}_{\L_t}}.$
\end{lemma}
We remark that the first of these estimates coincides with the
heuristic computation that since the collar around $\a$ has length
comparable to $\log[1/l(\a)]$, the derivative should be approximately
$-1/l(\a)$.

\medskip

We also need to bound the coefficient $\cosh u_j$ of $\mbox{d}d_j$ in
Equation~(\ref{eqn:cmswolpert}). Notice that the bound applies to all
the pants curves (not just the extremely short ones) in a short
marking.
  \begin{lemma}
\label{lem:ubound}
If the pants curves system $\cal P$ is short, then for all $j$,
$$|\cosh u_j | \mul
1.$$
\end{lemma}
\begin{proof} 
  Since ${\cal P}$ is short, by definition all curves $\a_j$ have
  length bounded above, and hence the length $d_j$ of the common
  perpendicular $\pi_j$ to $\tilde \a_j, \tilde \a_{j+1}$ is bounded
  below.

  First suppose that $\tilde \zeta$ intersects the infinite geodesic
  $\hat \pi_j$ that contains $\pi_j$.  In this case, $u_j = i \phi_j$,
  where $\phi_j$ is the angle between $\tilde \zeta$ and $\hat \pi_j$
  at their intersection point $o$.  Consider the case when $o$ is
  contained in the segment $\pi_j$. Let $x_j$ be the distance between
  $o$ and the endpoint $o_j$ of $\pi_j$ that lies on $\tilde \a_j$.
  Since $\tilde \zeta$ intersects both $\tilde \alpha_j$ and $\tilde
  \alpha_{j+1}$, the angle of parallelism formula gives
$$
|\tan \phi_j |< 1/\sinh x_j \quad\text{and}\quad |\tan \phi_j| <
1/\sinh({d_j - x_j}).
$$
Since at least one of $x_j $ and $ d_j - x_j$ is bounded below,
this gives a uniform upper bound on $ |\tan \phi_j| $. Thus $\phi_j$
is uniformly bounded away from $\pi/2$ and $ |\cosh u_j| = |\cos
\phi_j| $ is bounded below by a universal positive number.

Now consider the case when $o$ lies outside of $\pi_j$.  Let $o_j',
o'_{j+1}$ denote respectively, the points of intersection between
$\tilde \zeta_j$ and $\tilde \alpha_j, \tilde \alpha_{j+1}$ and let
$o_j$, $o_{j+1}$ denote respectively, the points of intersection
between $\pi_j$ and $\tilde \alpha_j, \tilde \alpha_{j+1}$. If
$d(o,o_j) \geq d(o,o_{j+1})$, then replace $\tilde \zeta$ with the
geodesic $\tilde \zeta'$ that passes through $o_j'$ and $o_{j+1}$ and
if $d(o, o_j) \leq d(o,o_{j+1})$, then replace $\tilde \zeta$ with the
geodesic $\tilde \zeta'$ that passes through $o_{j+1}'$ and $o_{j}$.
The angle $\phi'_j$ of intersection between $\tilde \zeta'$ and
$\pi_j$ satisfies $|\cos \phi_j| > |\cos \phi_j'|$.  We now run the
preceding argument with $\tilde \zeta'$ in place of $\tilde \zeta$ to
conclude $|\cos \phi_j'|$ is bounded below.

Now, suppose that $\tilde \zeta$ does not intersect $\hat \pi_j$.
Then $u_j$ is the hyperbolic distance from $\tilde \zeta$ to $\hat
\pi_j$. Denote by $p$ the point where the common perpendicular from
$\tilde \zeta$ to $\hat \pi_j$ meets $\hat \pi_j$; this point may lie
outside the segment $\pi_j$ between $\tilde \a_j, \tilde \a_{j+1}$.
Let $y_j, y_j'$ denote the (unsigned) distances from $p$ to $\tilde
\alpha_j, \tilde \alpha_{j+1}$ respectively. The quadrilateral formula
gives $\sinh y_j \sinh u_j = |\cos \theta_j|$ and $\sinh y_j' \sinh
u_j = |\cos \theta_{j+1}|$, where $\theta_j, \theta_{j+1}$ are the
angles between $\tilde \zeta$ and $\tilde \alpha_j, \tilde
\alpha_{j+1}$ respectively. Whether or not $ p \in \pi_j$, at least
one of $y_j $ and $y_j' $ is bounded below by $d_j/2$. Since there is
a uniform lower bound on $d_j$, it follows that $|\sinh u_j|$ and
hence $|\cosh u_j|$ is uniformly bounded above.  The result follows.
\end{proof}
  
\medskip

We now consider the second sum in Equation~(\ref{eqn:cmswolpert}). 

\begin{lemma} 
\label{lem:shift}
Let $\tilde \alpha_j$ be a lift of the curve $\a$ along which the
curve $\zeta$ has 
 shift coordinate $S_j=S_j (\zeta)$. Then 
$$ \frac{\dd S_j}{\dd l(\a)} = tw_{\L_t}(\zeta,\a) -
 s_{\a}(\L_t) \pm O(1),$$ where 
 $s_{\a}(\L_t)$ is the Fenchel-Nielsen twist along $\a$ at $\L_t$ as defined in
   Section~\ref{sec:fncoords}. 
\end{lemma} 
 \begin{proof} 
   Homotope the lift $\tilde \zeta$ to the piecewise geodesic path
   $\hat \zeta$ that runs along the successive lifts $\tilde \a_i$ and
   common perpendiculars $\pi_i$. The projection of $\hat
   \zeta$ to $S$ is homotopic to $\zeta$.  Then $S_j (\zeta)$ equals
   the signed distance that $\hat \zeta$ travels along $\tilde \a_j$.
    We need to express this shift   in a usable form.
  
    Recall the definition of Fenchel-Nielsen twist coordinates from
    Section~\ref{sec:fncoords}.  As above, we denote the pants curves
    in a short marking for $\L_t$ by ${\a_1}, \ldots,{\a_k} $.  The
    curve $\tilde \a_j$ forms the boundary of the lifts of two,
    possibly coincident, pairs of pants $P_j$ and $P_{j+1}$.  The
    projection $\hat \a_j$ of $\tilde \a_j$ to $P_j$ is bisected by
    the endpoints of the two seams of $P_j$ which join $\hat \a_j$ to
    each of the other two boundaries of $P_j$ (before identification
    in the surface $S$).  Likewise the projection $\hat \a_j'$ of
    $\tilde \a_j$ to $P_{j+1}$ is bisected by the endpoints of exactly
    two seams of $P_{j+1}$.
 
    The zero twist surface $\s_0= \s_0(l_{\a_1}, \ldots,l_{\a_k} ) $
    is formed by gluing $P_j$ to $P_{j+1}$ along $\hat \a_j$ and $\hat
    \a_j'$ in such a way as to match these two pairs of points. Thus
    on $\s_0$, the distance along $\tilde \a_j$ between incoming and
    outgoing perpendiculars $\pi_j$ and $\pi_{j+1}$ may be expressed
    in the form $ n_j(\zeta) l(\a)/2 + e_j(\zeta) $, where $n_j
    (\zeta) \in \mathbb{Z}$ is the (signed) number of seams $\hat
    \zeta$ intersects along $\tilde \a_j$ and $e_j(\zeta)$ is an error
    term which allows for the possibility that $\pi_j, \pi_{j+1}$ may
    not be seams of $P_j$ and $P_{j+1}$, but rather common
    perpendiculars from $\hat \a_j$ or $\hat \a_j'$ to itself.  In all
    cases however, $|e_j(\zeta)| < l_\s(\a_j)$ and $e_j(\zeta)$
    depends only on the geometry of $P_j$ and $P_{j+1}$,
    see~\cite{cmswolpert}.

    Now at $\L_t$, the incoming and outgoing perpendiculars $\pi_j$
    and $\pi_{j+1}$ are further offset by
    $\tau_{\a_j}(\L_t)=l_{\L_t}(\a)s_{\a_j}(\L_t)$ giving the formula
   $$S_j(\zeta) = \frac{1}{2}n_j(\zeta) l_{\L_t}(\a) + e_j(\zeta) 
+ \tau_{\a_j}(\L_t),$$
   see also Section 4.2 of~\cite{cmswolpert}.
      
   We can now proceed to estimate $\dd S_j/\dd l(\a)$. Since the
   partial derivatives are taken with respect to the coordinates
   $(l_\s(\a),\tau_{\a}(\s))$, the term $\dd \tau_{\a_j}/\dd l(\a)$
   vanishes.  To avoid an unpleasant calculation, we get rid of the
   term $e_j(\zeta)$ as follows.  Modify $\hat \zeta$ to a path which
   still runs along the lifts of the pants curves and their common
   perpendiculars, but which never goes along a perpendicular from a
   lift of $\a$ to itself.  Specifically, let $\pi_j$ be such a common
   perpendicular which projects to a pair of pants $P$ one of whose
   boundary components is $\a$. Let $\beta$ be one of the other
   boundary components and let $\eta$ be the perpendicular from $\a$
   to $\beta$.  The projection of $\pi_j$ to $P$ is homotopic, with
   fixed endpoints, to an arc which runs along $\a$, then along
   $\eta$, then along $\beta$, back along $\eta$, finally back to the
   final point on $\a$, see Figure~\ref{fig:pants}.
\begin{figure}[htb]
 \centerline{\epsfbox{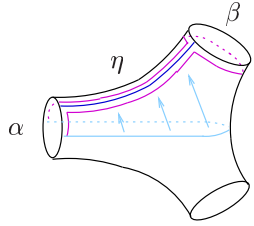}}
\caption{Homotoping the perpendicular.}
\label{fig:pants}
\end{figure}
Modify $\hat \zeta$ by replacing $\pi_j$ by the lift of this alternate
path.  Doing this in each instance gives a replacement for $\hat
\zeta$, with respect to which one can define all quantities occurring
in (\ref{eqn:cmswolpert}) as before. The derivation of
Equation~(\ref{eqn:cmswolpert}) in~\cite{cmswolpert} will still work
for this new path.  Denoting the newly defined shift also by $S_j$, we
thus have $\dd S_j/\dd l(\a) = n_j/2 \pm 1$.

We claim that $n_j/2 = tw_{\s_0}(\zeta,\a) \pm O(1)$.  By definition,
$\hat \zeta$ traverses lifts of the pants curves to $\HH$ in the same
order as $\tilde \zeta$.  Thus the segment of $\hat \zeta$ running
along the lift $\tilde \a = \tilde \a_j$ to $\HH$ is the interval
between the footpoints $Q_j,Q_{j+1}$ of the perpendiculars $\pi_j$ and
$\pi_{j+1}$ from $\tilde \a_{j-1}$ and $\tilde \a_{j+1}$ (the lifts of
pants curves adjacent to $\tilde \a_j$) to $\tilde \a_j$.  Now $Q_j$
lies within the interval on $\tilde \a_j$ bounded by the footpoints of
the perpendiculars from the two endpoints of $\tilde \a_{j-1}$ on
$\partial \HH$ to $\tilde \a_j$; and similarly for $ Q_{j+1}$. Thus
our claim follows as in~\cite{minskyproduct} Lemma 3.1, see the
discussion in Section~\ref{sec:twistshyperbolic}.  The proof of the
present lemma can now be completed by applying
Lemma~\ref{lem:minskytwist}.
\end{proof}

We can now put the above results together to obtain an estimate of
$\partial l(\nu)/\partial l(\a)$.  Let $\{ \pi_j \}_{j \in J}$ be the
subset of perpendiculars whose projections are contained in a common
pair of pants with $\a$ but are disjoint from $\a$. Then by
Lemmas~\ref{lem:hexagon*}(ii) and ~\ref{lem:ubound}, we have
$$\sum_{j \in J} \cosh u_j \frac{\partial d_j}{\partial l(\a)} \mul
\sum_{j \in J} \frac{\partial d_j}{\partial l(\a)} \mul \sum_{j \in J}
l(\a).$$
In the case that we have a measured lamination $\nu$ instead
of a curve $\zeta$, we obtain by the same reasoning
$$
\sum_{j \in J} \cosh u_j \frac{\partial d_j}{\partial l(\a)} \mul
\sum_{j \in J} l(\a) \cdot w_\nu(\pi_j),$$
where $w_\nu(\pi_j)$ is the
$\nu$-weight of the leaves of $\tilde \nu$ that go from $\tilde \a_j$
to $\tilde \a_{j+1}$.  Let us denote
\begin{equation}
\label{eqn:defofDel}
\Delta_\nu(\a) = \sum_{j \in J} w_\nu(\pi_j).
\end{equation}
By applying Equation~(\ref{eqn:twistangle}) and Lemmas~\ref{lem:hexagon*}
--~\ref{lem:shift} to Equation~(\ref{eqn:cmswolpert}), we obtain the
following:
\begin{lemma}
\label{lem:lengthderiv} Let $\a$ be a curve in a short pants
decomposition $\cal P_{\L_t}$ of $\L_t$  and let $\nu$ be a measured
lamination transverse to $\a$, with average intersection angle
$\Theta$. If $\a$ is extremely short and neither horizontal nor
vertical, then using coordinates $(l (\a_i), \tau_{\a_i}) $ relative
to ${\cal P}_{\L_t}$, we have $\partial l(\nu)/\partial l(\a) = -A + B + C$,
where
\begin{align*}
  A &\mul i(\nu,\a)\frac{1}{l(\a)}, \qquad B \mul \Delta_\nu(\a)
  l(\a) \qquad\text{and} \\
  C &\mul i(\nu,\a)\left( (tw_{\L_t}(\nu,\a) - s_\a(\L_t)) \cos \Theta
    + O(1) \right).
\end{align*}
\end{lemma}
  
\subsection{Proof of the main estimates} 
We are ready to prove our main results Theorems~\ref{thm:shortonL}
and~\ref{thm:twistonL}.  If we apply Equation~(\ref{eqn:1form}) to
$\partial/\partial l(\a)$, then by Lemma~\ref{lem:lengthderiv}, we
obtain
\begin{equation}
\label{eqn:abcis0}
0=\frac{\dd l(\nu^+_t)}{\dd l(\a)} 
+   \frac{\dd l(\nu^-_t)}{\dd l(\a)}  = -(A^+ + A^-)
+ B^+ + B^- + C^+ + C^-, 
\end{equation}
where
\begin{align*}
A^\pm &\mul \frac{i(\nu_t^\pm,\a)}{l(\a)},\qquad
B^\pm \mul \Delta_{\nu^\pm_t}(\a)l(\a) \qquad\text{and}\\
C^\pm &\mul i(\nu_t^\pm,\a) \left( (tw_{\L_t}(\nu^\pm,\a)  
- s_\a(\L_t) )\cos \Theta^{\pm}_t + O(1) \right).
\end{align*}
Since $A^+ + A^-= B^+ + B^- + C^+ + C^-$, we get
\begin{equation}
\label{pre}
\frac{1}{l(\a)} \mul
\frac{\Delta_{\nu^+_t}(\a) +
\Delta_{\nu^-_t}(\a) }{i(\nu^+_t,\a) + i(\nu^-_t,\a)}l(\a)
+ \frac{ C^+ + C^-}{i(\nu^+_t,\a) + i(\nu^-_t,\a)}.
\end{equation}
Notice that the term $C^+ + C^-$ simplifies: 
defining  
$$D^{\pm} = C^{\pm} + i(\nu_t^\pm,\a) s_\a(\L_t)\cos \Theta^{\pm}_t,$$
it follows immediately from Equation~(\ref{eqn:kerckhoff}) that $C^+ +
C^- = D^+ + D^-$.

\begin{lemma} 
\label{lem:3} 
Let $H_t(\a)$ be
  defined as in Equation~(\ref{eqn:defofH}). Then we have
  $$H_t(\a) \asymp \frac{\Delta_{\nu^+_t}(\a) + \Delta_{\nu^-_t}(\a)
  }{i(\nu^+_t,\a) + i(\nu^-_t,\a)}.$$
\end{lemma}
\begin{proof}
  The strands of $\nu_t^\pm$ which intersect pants adjacent to $\a$
  but which are disjoint from $\a$, must intersect one of the curves
  in $\cal B$. Hence, by definition of $\Delta_{\nu^\pm_t}(\a)$, we
  have
  $$
  \frac{\Delta_{\nu^+_t}(\a) + \Delta_{\nu^-_t}(\a) }{i(\nu^+_t,\a)
    + i(\nu^-_t,\a)} \prec \sum_{\beta \in \cal B}
  \frac{l_{q_t}(\beta)}{l_{q_t}(\a)} \mul H_t(\a).$$
  To prove the
  inequality in the other direction, let $\beta \in \cal B$ be the
  curve that realizes the maximum in the definition of $H_t(\a)$.  For
  $\nu = \nu^\pm_t$, let $\nu_{\b \a}, \nu_{\b \b}, \nu_{\b \gamma}$
  be the collections of strands of $\nu$ that run between $\beta$ and
  $\a$, from $\beta$ to itself, and between $\beta$ and $\gamma$,
  respectively.  (As usual, there are different possible
  configurations of strands in each pants, in particular $\nu_{\b \b}$
  may be empty. The inequalities which follow are however valid in all
  cases.)
  Denote the $\nu$-weight of these by $w(\nu_{\b \a})$, $w(\nu_{\b
    \b})$, and $w(\nu_{\b \gamma})$, respectively.  Then
\begin{align*}
  H_t(\a) & = \frac{l_{q_t}(\beta)}{l_{q_t}(\a)} \mul
  \sum_{I=+,-}\frac{w(\nu^I_{\b \a}) + w(\nu^I_{\b \b})
    + w(\nu^I_{\b\gamma})}{l_{q_t}(\a)} \\
  & \prec \sum_{I=+,-}\frac{w(\nu^I_{\b \b}) + w(\nu^I_{\b
      \gamma})}{l_{q_t}(\a)} \prec \frac{\Delta_{\nu^+_t}(\a) +
    \Delta_{\nu^-_t}(\a) }{i(\nu^+_t,\a) + i(\nu^-_t,\a)}. \qedhere
\end{align*}
\end{proof}

\begin{proof}[Proof of Theorem~\ref{thm:shortonL}]
 Lemma~\ref{lem:3}, Equation~(\ref{pre}), and the remark following
  gives $
  1/l_{\L_t}(\a) \asymp G_t +H_t l_{\L_t}(\a), $ where $H_t=H_t(\a)$
  and
$$G_t= G_t(\a)= \frac{D^+ + D^-}{i(\nu^+_t,\a) + i(\nu^-_t,\a)}. $$
Hence, we must have either
$$
\frac{1}{l_{\L_t}(\a)} \asymp
 G_t \hspace{0.3cm} \mbox{ or }\hspace{0.3cm}
 \frac{1}{l_{\L_t}(\a)} \asymp  {H_t}l_{\L_t}(\a),
$$
from which we obtain
$$
\frac{1}{l_{\L_t}(\a)} \asymp \max \{ G_t,
\sqrt{H_t} \}.
$$

We simplify the expression for $G_t$ as follows.  By the discussion
following Equation~(\ref{eqn:twistangle}) we see that either $
tw_{\L_t}(\nu^{+},\a) \cos \Theta^{+}_t$ is positive, or
$|tw_{\L_t}(\nu^{+},\a) \cos \Theta^{+}_t| = O(1)$, and likewise for
$\nu^-$. Also note that $tw_{\L_t}(\nu^{+},\a)$ and
$tw_{\L_t}(\nm,\a)$ are either $O(1)$ or have opposite signs, so that
$$\da = Tw_{\L_t}(\nu^+,\a) + Tw_{\L_t}(\nm,\a) \pm O(1).$$ 
As before, let $T = t-t_\a$.  Then by applying Equation~(\ref{eqn:kerckhoff}), 
we get
\begin{align}
  G_t(\a) &= \frac{e^T \left| \cos \Theta^+_t \right|
    Tw_{\L_t}(\nu^+,\a) + e^{-T} \left| \cos \Theta^-_t
    \right|Tw_{\L_t}(\nm,\a)}
  {e^T + e^{-T}} + O(1)  \nonumber\\
&\asymp 
\frac{e^{T} \left| \cos \Theta^+_t \right| \da}{e^{T} + e^{-T}}  =  
\frac{e^{-T} \left|\cos \Theta^-_t \right| \da}{e^{T} + e^{-T}}.
\label{G}
\end{align}

This almost completes the proof, except it remains to be shown that if
$1/l_{\L_t}(\a) \asymp G_t$, then 
$$G_t \asymp e^{-2|t-t_\a|}\da.$$
By Equation~(\ref{G}), it is
sufficient to show that there is some constant $c>0$, independent of
$\a$, such that $\left| \cos \Theta^-_t \right| > c$ whenever $T>0$
and $\left| \cos \Theta^+_t \right|> c$ whenever $T<0$.

Our assumption that $1/l_{\L_t}(\a) \asymp G_t$ and the fact that
$l_{\L_t}(\a)$ is sufficiently small, together with Equation~(\ref{G})
imply that
$$\frac{1}{l_{\L_t}(\a)}\stackrel{{}_\ast}{\prec} Tw_{\L_t}(\nm,\a) +
Tw_{\L_t}(\np,\a).$$
Let $X_t = Tw_{\L_t}(\nm,\a)l_{\L_t}(\a)$ and
$Y_t = Tw_{\L_t}(\np,\a)l_{\L_t}(\a)$. The above inequality states
that
\begin{equation}
\label{XplusY}
 X_t + Y_t \stackrel{{}_\ast}{\succ} 1.
\end{equation}
If $T>0$, then by Equation~(\ref{eqn:kerckhoff}), $\left| \cos
  \Theta^{-}_t \right| > \left| \cos \Theta^{+}_t \right|$ so $X_t >
Y_t - O(l_{\L_t}(\a))$ by Equation~(\ref{eqn:twistangle}).  Thus,
reducing the value of the upper bound $\ep_0$ on $l_{\L_t}(\a)$ if
necessary, it follows from Equation~(\ref{XplusY}) that $X_t$ is
bounded below by some positive constant, and thus the same is true of
$|\cos{\Theta^-_t}|$.  The analogous statement holds for $|\cos
\Theta^+_t|$ when $T<0$.
\end{proof}

\begin{proof}[Proof of Theorem~\ref{thm:twistonL}]
  From Equation~(\ref{eqn:kerckhoff}), $|\cos \Theta^\pm_t| \leq
  e^{-2|T|}$.  It follows from Equation~(\ref{eqn:twistangle}) that if
  $T \gg 0$ then
  $$Tw_{\L_t}(\np,\a) \, l_{\L_t}(\a) \prec e^{-2T}.$$
  The argument
  for $T \ll 0$ is similar.  Now suppose that $|T|= O(1)$. Since
  $$
  Tw_{\L_t}(\nu^\pm,\a) \, l_{\L_t}(\a) \prec \da \, l_{\L_t}(\a)$$
  and since by Theorem~\ref{thm:shortonL},
  $$\da \, l_{\L_t}(\a) \prec e^{2|T|},$$
  the result follows.
\end{proof}

\subsection{Estimates for horizontal and vertical short curves on $\L$}
\label{sec:vertonL}
As in Section~\ref{sec:vertonG}, we need the analogue of
Theorems~\ref{thm:shortonL} and \ref{thm:twistonL} for extremely short
curves $\a$ which are either horizontal or vertical.  As in that
section, assume $\a$ is vertical so that $i(\a,\nu^-) = 0$.

As before, we shall obtain the estimates by applying
Equation~(\ref{eqn:1form}) to $\frac{\dd}{\dd \tau_\a}, \frac{\dd}{\dd
  l(\a)}$.  Since $\a$ is vertical,
$$0=\frac{\dd l(\nu_t^+)}{\dd \tau_\a}+\frac{\dd l(\nu_t^-)}{\dd
  \tau_\a} =\frac{\dd l(\nu_t^+)}{\dd \tau_\a} = \int_\a \cos \theta^+
d\nu_t^+.$$
Hence Equation~(\ref{eqn:kerckhoff}) is replaced by $\cos
\Theta^+_t = 0$.  Furthermore, Equation~(\ref{eqn:twistangle}) gives
$|tw_{\L_t}(\nu^+,\a)| = O(1)$.

Let $m^-(\a)$ be the weight on $\a$ of $\nu^- =\nu_0^-$, in other
words, $\nu^- =m^-(\a) \a + \eta$, where $\eta$ has support disjoint
from $\a$.  Then following the line of discussion in
Section~\ref{sec:differentiation}, it is easy to check that
$$
\frac{\dd l(\nu_t^-)}{\dd l(\a)} = e^{-t} m^-(\a) +
\Delta_{\nu_t^-}(\a) l(\a).$$
Hence, in place of Equation~(\ref{eqn:abcis0}), we obtain
$$0=\frac{\dd l(\nu_t^+)}{\dd l(\a)}+\frac{\dd l(\nu_t^-)}{\dd l(\a)}
= -A^+ + B^+ + C^+ + e^{-t} m^-(\a) + B^-,$$
where $A^+, B^\pm, C^+$ are defined as before. Since
$$
C^+ \mul i(\nu^+_t,\a) \left[ \big( tw_{\L_t}(\nu,\a) -
s_\a(\L_t) \big) \cos \Theta_t^+ + O(1) \right] \mul i(\nu^+_t,\a),
$$ 
we get
$$\frac{1}{l_{\L_t}(\a)} \asymp H_t(\a) l_{\L_t}(\a) + 
\frac{e^{-t} m^-(\a)}{i(\nu_t^+,\a)}=H_t(\a) l_{\L_t}(\a)+
e^{-2t}\frac{m^-(\a)}{i(\nu^+,\a)}.$$
Thus we obtain:
\begin{theorem}
 \label{thm:shortonL1}
 Let $\a$ be a curve which is vertical on $S$. If $\a$ is extremely
 short in $\L_t$, then
$$\frac{1}{l_{\L_t}(\a)} {\asymp}  
\max \Big\{ e^{-2t}\frac{m^-(\a)}{i(\nu^+,\a)}, \sqrt{H_t(\a)} \Big\}. 
$$
If $\a$ is horizontal, the estimate is the same except that the
first term is replaced by $e^{2t}m^+(\a)/i(\nu^-,\a)$, where now
$m^+(\a)$ is the weight on $\a$ of $\nu^+$.
\end{theorem}

\begin{theorem}
 \label{thm:twistonL1} 
 If $\a$ is extremely short in $\L_t$, then the twist satisfies
 $Tw_{\L_t}(\nu^+,\a) = O(1)$ if $\a$ is vertical and
 $Tw_{\L_t}(\nu^-,\a) = O(1)$ if $\a$ is horizontal.
\end{theorem}

\section{Comparing $\L_t$ and $\G_t$}\label{sec:finalcompare}

In this section we prove our final results. We compare the geometry of
$\L_t$ and $\G_t$ by looking at their respective thick-thin
decompositions.  Specifically, we prove that
\begin{equation}
\label{HisK}
H_t(\a) \asymp K_t(\a).
\end{equation} Combined
with Theorems~\ref{prop:shortingeo} and~\ref{thm:shortonL}, this
completes the proof of Theorems~\ref{shortonL} and \ref{shortonboth}.
We show further in Theorem~\ref{thm:markingsclose} that on
corresponding thick components, the two metrics $\L_t$ and $\G_t$
almost coincide.  Combining this with the information about twisting
given in Theorems~\ref{prop:quadtwistest} and \ref{thm:twistonL}, we
can then estimate the \Teich distance between $\L_t$ and $\G_t$, thus
completing the proof of Theorem~\ref{distance}
(Theorem~\ref{thm:finalcomparison}).

As explained in the Introduction, the logical flow in the proof of
Equation~(\ref{HisK}) is not straightforward. We first show relatively
easily in Proposition~\ref{cor:shorteronL} that $H_t(\a) \succ
K_t(\a)$.  The key point in proving the other half of
Equation~(\ref{HisK}) is Proposition~\ref{lem:minimalonQ}, which shows
that the metric $\L_t$ not only minimizes $l_{\s}(\nu^+_t ) +
l_{\s}(\nu^-_t)$, but that it also in a suitable coarse sense
minimizes the contribution to the sum made by the parts of
$\nu^{\pm}_t$ which lie in the thick part of $\L_t$.  This is proved
in Section~\ref{sec:subsurfaces}.  In Section~\ref{sec:correspondence}
we use Proposition~\ref{lem:minimalonQ} to deduce that a curve that is
extremely short in $\G_t$ is also extremely short in $\L_t$
(Proposition~\ref{bddbelow}) .  This is used in proving
Theorem~\ref{thm:markingsclose} mentioned above, from which in
Section~\ref{sec:compshtlength} we are finally able to show that
$H_t(\a) \prec K_t(\a)$ (Proposition~\ref{lem:HandK}).

\subsection{Curves are shorter in $\L(\np,\nm)$}
\begin{proposition}
\label{cor:shorteronL} 
If $\a$ is extremely short in $\mathcal G_t$, then
$H_t(\a) \succ K_{t}(\a)$. Therefore,
$$\frac{1}{l_{\L_t}(\a)} \succ \frac{1}{l_{\mathcal G_t}(\a)}.$$
\end{proposition}
\begin{proof}
  Once we show that $H_t(\a) \succ K_{t}(\a)$, the second statement
  follows from Theorems~\ref{prop:shortingeo} and \ref{thm:shortonL}.
  
  The only case of interest is when $ K_{t}(\a)$ is large. Let
  $E_t(\a)$ be one of the expanding annuli $E_i(\a)$ around $\a$ of
  larger modulus, defined as in the discussion preceding
  Corollary~\ref{cor:lengthest}. Denote the inner and outer boundary
  curves of $E_t(\a)$ by $\dd_0$ and $\dd_1$.  Let $\omega$ be an
  essential arc from $\a$ to itself such that $l_{q_t}(\omega) = 2
  d_{q_t}(\dd_0,\dd_1)$, where as usual $d_{q_t}$ denotes distance in
  the $q_t$-metric. The annulus $E_t(\a)$ intersects the
  $q_t$-representative $\hat Y$ of a thick component $Y$ of
  $(S,\mathcal G_t)$ adjacent to $\a$.  Let us first suppose that
  $\omega$ is contained in $\hat Y$. A small regular neighborhood of
  $\a \cup \omega$ has boundary consisting of $\a$ and two curves,
  $\zeta_1, \zeta_2$, which together with $\a$ bound a pair of pants.
  Therefore, either both $\zeta_1$ and $\zeta_2$ are contained in
  ${\cal B}$ or one of these two curves must intersect a curve in
  ${\cal B}$ transversely.  (As in Equation~(\ref{eqn:defofH}), ${\cal
    B}$ is the set of pants curves in a short pants decomposition
  ${\cal P}_{\L_t}$ which are boundaries of pants adjacent to $\a$.)

First consider the case when $\zeta_1,\zeta_2 \in {\cal B}$. If
either $\zeta_1$ or $\zeta_2$ is non-peripheral in $Y$, then
by definition of $H_t(\a)$ and $\lambda_Y$, 
$$
H_t(\a) =
\max_{\beta \in \cal B} \left\{\frac{l_{q_t}(\beta)}{l_{q_t}(\a)} \right\} \geq 
\max \left\{\frac{l_{q_t}(\zeta_1)}{l_{q_t}(\a)}, 
\frac{l_{q_t}(\zeta_2)}{l_{q_t}(\a)} \right\} \geq 
\frac{\lambda_Y}{l_{q_t}(\a)}.
$$
If both $\zeta_1,\zeta_2$
are peripheral in $Y$, then
\begin{equation}
\label{peripheral}
H_t(\a) \succ
\frac{l_{q_t}(\zeta_1)+l_{q_t}(\zeta_2) + l_{q_t}(\a)}{l_{q_t}(\a)}
\asymp \frac{\lambda_Y}{l_{q_t}(\a)}.
\end{equation}

Now consider the case when either $\zeta_1$ or $\zeta_2$ intersects a
curve $\beta \in {\cal B}$ transversely. Note that for $i=1,2$,
$$l_{q_t}(\zeta_i) \leq l_{q_t}(\omega) + l_{q_t}(\a) \leq 4\,
\mbox{diam}_{q_t}(\hat Y) \mul \lambda_Y.$$
Suppose that $\beta$
intersects $\zeta_i$ and that $\zeta_i$ is not peripheral in $Y$. Then
it follows from the above inequality and the definition of $\lambda_Y$
that $l_{q_t}(\zeta_i) \mul \lambda_Y$.  Since by
Theorem~\ref{thm:scalefactor} we have $l_{q_t}(\zeta_i) \mul \lambda_Y
l_{\mathcal G_t}(\zeta_i)$, we see that $l_{\mathcal G_t}(\zeta_i)
\mul 1$.  Then, by the collar lemma for quadratic differentials
\cite{rafi1}, we have
$$l_{q_t}(\beta) \stackrel{{}_\ast}{\succ} l_{q_t}(\zeta_i) \mul
\lambda_Y.$$
The only remaining possibility to consider is when both
$\zeta_1,\zeta_2$ are peripheral so that $Y$ is a pair of pants.  If
$\beta$ intersects $\zeta_i$, then since $\beta$ must pass through an
annulus around $\zeta_i$ which has large modulus, we conclude
$l_{q_t}(\beta) \succ l_{q_t}(\zeta_i)$. If $\beta$ does not intersect
$\zeta_i$, every arc $\eta$ of $\beta \cap \hat Y$ has both endpoints
on the other curve $\zeta_j$, $j\neq i$.  The endpoints of $\eta$
divide $\zeta_j$ into two arcs, one of which together with $\eta$
forms a curve homotopic to $\zeta_i$.  Since $\beta$ passes through an
annulus of large modulus around $\zeta_j$, this implies
$$l_{q_t}(\zeta_i) \leq l_{q_t}(\eta) + l_{q_t}(\zeta_j)
\stackrel{{}_\ast}{\prec} l_{q_t}(\beta) + l_{q_t}(\beta).$$
Either
way, we have $l_{q_t}(\beta) \stackrel{{}_\ast}{\succ}
l_{q_t}(\zeta_i)$ and thus $l_{q_t}(\beta) \stackrel{{}_\ast}{\succ}
l_{q_t}(\zeta_1) + l_{q_t}(\zeta_2)$.

Finally, if the original arc $\omega$ was not contained in $\hat Y$,
we can replace it with an arc that is contained in $\hat Y$ of
comparable length as in the proof of Theorem~\ref{expanding}, and run
the same argument.
\end{proof}

\subsection{Length estimates on subsurfaces}
\label{sec:subsurfaces}

The object of this section is to prove
Proposition~\ref{lem:minimalonQ}.  We begin with estimates that are
necessary to analyze the contribution to the length of a lamination
associated to the thick part of the surface $S$.  Thus if $(S,\s)$ is
a hyperbolic surface and $Y \subset S$ is a subsurface of the thick
part, we want to find an approximation to $l_{\s}(\nu^{\pm} \cap Y)$.
To consider the problem in general, we consider $l_\s(\nu^\pm \cap Q)$
for a subsurface $Q$ with geodesic boundary.  Suppose that $\zeta$ is
a geodesic that intersects $Q$ but is not entirely contained in $Q$.
The essential idea is that we can approximate $\zeta \cap Q$ by
piecewise geodesic arcs homotopic to $\zeta \cap Q$, which alternately
run along arcs perpendicular to $\partial Q$ and parallel to $\partial
Q$.  The length of the parallel portion is determined by the twisting
of $\zeta$ about the curves in $\partial Q$, while the portion
$\zeta_Q$ perpendicular to $\partial Q$ is defined and estimated as
explained below.

Let $\ua$ be a collection of disjoint simple closed geodesics on
$(S,\s)$ and let $Q$ be a totally geodesic surface which is the metric
completion of a component of $S\setminus \ua$.  (It is possible for
two distinct boundary components of $Q$ to be identified in $S$ to a
single curve $\a \in \ua$, so strictly speaking, $Q$ is not a
subsurface of $S$.)  If $\eta$ is an essential geodesic arc with
endpoints on $\partial Q$, let $\eta_Q$ be the shortest arc in $Q$
that is freely homotopic to $\eta$, relative to $\dd Q$. In this case,
clearly $\eta_Q$ is orthogonal to $\dd Q$. If $\varphi$ is a measured
geodesic lamination whose support is entirely contained in $Q$, let
$\varphi_Q = \varphi$. For convenience we allow the possibility that
the support of $\varphi$ contains components of $\dd Q$, remarking
that this is not quite the same as the definition
in~\cite{minskyproduct}.

Suppose $\xi$ is a measured geodesic lamination on $S$.  Then the
intersection $\xi \cap Q$ is a union of components of $\xi$ that are
entirely contained in $Q$ and arcs with both endpoints on $\dd Q$.  If
$\eta$ is an arc of $\xi \cap Q$, let $n(\eta_Q) $ denote the
transverse measure of arcs in the homotopy class $[\eta_Q]$. The
\emph{orthogonal projection} of $\xi$ into $Q$ is $ \xi_Q= \sum
n(\eta_Q)\eta_Q + \sum \varphi_Q$, where the first sum is taken over a
representative $\eta_Q$ from each class of arcs in $\xi \cap Q$ and
the second sum is taken over all components $\varphi$ of $\xi$ that
are entirely contained in $Q$. Define
$$l_{\s}(\xi_Q) = \sum n(\eta_Q)l_{\s}(\eta_Q) + \sum l_\s(\varphi_Q).$$

If all curves in $\partial Q$ are of uniformly bounded length, then we
have the following estimate of $l_\s(\xi \cap Q)$ in terms of $\xi_Q$
and $Tw_{\sigma}(\xi,\a)$:
\begin{lemma}
\label{lem:estimate}
Suppose $l_\s(\a_j) <\ell$ for every component $\a_j$ of $\dd Q$. Then
there exists a constant $K=K(\ell)$ such that for any measured
lamination $\xi$ on $S$:
$$
\left|\,l_\s(\xi \cap Q) -\bigg[ l_\s(\xi_Q) +
\sum  l_\s(\a_j) \frac{ Tw_\s(\xi,\a_j)}{2} 
i(\xi, \a_j)\bigg]\,\right| \leq K \, i(\xi, \dd Q),
$$
where the sum is taken over all $\a_j$ that intersect $\xi$
transversely.
\end{lemma}
\noindent For a proof, see the Appendix.
The next lemma can be proved similarly, applying the same property of
hyperbolic triangles. We omit the proof.
\begin{lemma}
\label{lem:thinlength} 
Suppose $l_\s(\a) < \ep_{0}$. Let $A$ be an embedded annulus in
$(S,\s)$ such that one component of $\dd A$ is the geodesic $\a$ and
the other a hyperbolically equidistant curve of length $\ep_0$.  Then
there is a uniform constant $K$ such that for any measured lamination
$\xi$ on $S$ that intersects $\a$ transversely:
$$
\left|\,l_\s(\xi \cap A) - \left[ \log \frac{1}{l_\s(\a)} 
+ l_\s(\a)\frac{Tw_\s(\xi,\a)}{2} \right] i(\xi, \a)\,\right|
\leq K \, i(\xi, \a).
$$
\end{lemma}
Here, $ \log [1/l_\s(\a)]$ approximates the width
of $A$, up to a bounded additive error. \\

We will now apply Lemmas~\ref{lem:estimate} and \ref{lem:thinlength}
to prove Proposition~\ref{lem:minimalonQ} below, which in turn is the
key step to proving Theorem~\ref{thm:markingsclose}.

For $\rho > 0$ and a hyperbolic metric $\s \in \ts$, define
$${\cal S}_\rho(\s) = \{ \a \in \cal S : l_\s(\a) < \rho \}.$$
Proposition~\ref{cor:shorteronL} implies that if $\a$ is extremely
short in $\mathcal G_t$, then we can choose a constant $\ep < \ep_0$,
depending only on $\ep_0$, such that if $l_{\mathcal G_t}(\a) < \ep$,
then $l_{\L_t}(\a) < \ep_0$.  In other words, we can choose $\ep <
\ep_0$ so that

\begin{equation}
\label{star}
{\cal S}_\ep(\mathcal G_t) \subset {\cal S}_{\ep_0}(\L_t).
\end{equation}

Now, let $Q=Q_t$ be a component of $S \setminus {\cal S}_\ep(\mathcal
G_t)$.  The metric $\L_t$ naturally endows $Q$ with the structure of
hyperbolic surface with geodesic boundary, which by the above,
satisfies $l_{\L_t}(\a) < \ep_0$ for all components $\a$ of $\dd Q$.
Henceforth, fix a constant $c$ that satisfies $\ep_0 < c < \ep_\M$.
For $\s = \mathcal G_t, \L_t$, or in general, any metric that makes
$Q$ a hyperbolic surface with geodesic boundary components that are
extremely short, define $C(\a,\s)$ to be the collar of $\a$ in
$(Q,\s)$ such that one component of $\dd C(\a)$ is (the geodesic
representative of) $\a$, and the other, the equidistant curve of
length $c$.  Because $c < \ep_\M$, the collars are all disjoint from
one another. Let $(Q_T,\s)$ be the metric subsurface of $Q$ defined
by:
$$(Q_T,\s) = (Q,\s) \setminus \coprod_{\a \in \dd Q} C(\a,\s).$$
In particular, every component of $\dd Q_T$ has length $c$.

Since $\L_t$ is on the line of minima, we have 
$$l_{\L_t}(\nu^+_t ) + l_{\L_t}(\nu^-_t) \leq l_{\mathcal G_t}(\nu^+_t ) +
l_{\mathcal G_t}(\nu^-_t).$$
The contribution to this inequality from $Q_T$ is
given as follows:
\begin{proposition}
\label{lem:minimalonQ}
If $Q$ is a component of $S \setminus {\cal S}_\ep(\mathcal G_t)$,
then
$$l_{\L_t}(\nu^+ \cap Q_T ) + l_{\L_t}(\nu^- \cap Q_T)
\stackrel{{}_\ast}\prec l_{\mathcal G_t}(\nu^+ \cap Q_T) + 
l_{\mathcal G_t}(\nu^- \cap Q_T).$$
\end{proposition}
 
To prove Proposition~\ref{lem:minimalonQ}, we need the following two
lemmas.
\begin{lemma}
\label{thesame}
Suppose $(Q,\s)$ is a hyperbolic surface with geodesic boundary such
that $l_\s(\a) < \ep_0$ for all $\a \in \dd Q$. Then for any measured
geodesic lamination $\xi$ on $(S,\s)$, we have
$$l_\s(\xi \cap Q_T) \mul l_\s(\xi_Q \cap Q_T).$$
\end{lemma}
\begin{proof}
  We consider $l_\s(\xi \cap Q_T)=l_\s(\xi \cap Q) - l_\s(\xi \cap
  \coprod C(\a))$ and apply
  Lemmas~\ref{lem:estimate},~\ref{lem:thinlength} to the right-hand
  side. The proof is completed by observing that
$$
l_\s(\xi_Q \cap Q_T) = l_\s(\xi_Q ) - \sum_{\a \in \dd Q} 
\log \frac{1}{l_\s(\a)} i(\xi, \a) + O\big( i(\xi,\dd Q) \big)
$$
and that $l_\s(\xi \cap Q_T) \stackrel{{}_\ast}{\succ} i(\xi,\dd
Q_T) = i(\xi,\dd Q)$, due to the fact that every component of $\dd
Q_T$ has an annular neighborhood of definite width.
\end{proof}

\begin{lemma}
\label{bilipschitz}
Suppose that $(Q,\s)$ and $(Q,\s')$ are two hyperbolic surfaces with
geodesic boundary whose boundary components are all extremely short.
Suppose that there is a short pants decomposition of $(Q,\s)$ with
respect to which the Fenchel-Nielsen coordinates for $(Q,\s)$ and
$(Q,\s')$ agree, except possibly for the lengths and twists
corresponding to components of $\dd Q$. Then for any simple closed
curve or arc $\eta$ with endpoints on $\dd Q$,
$$l_\s(\eta_Q \cap Q_T) \mul l_{\s'}(\eta_Q \cap Q_T).$$
\end{lemma}
\begin{proof}
  This is essentially the same as a discussion in
  Minsky~\cite{minskyproduct} page 283.  The idea is that there is a
  $K$-bilipschitz homeomorphism $(Q_T,\s) \to (Q_T,\s')$ with constant
  $K$ depending only on $\ep_0$. To see this, cut $Q$ along the pants
  curves into pairs of pants and further cut each pair of pants into
  hexagons. Corresponding hexagons in the two surfaces have the same
  side lengths, except those whose edges form part of $\partial Q$.
  Now truncate those hexagons which have an edge on $\dd Q$ by cutting
  off the collar round $\dd Q$ in such a way that the boundary of the
  truncated hexagon is the corresponding component of $\dd Q_T$. By
  our construction, the non-geodesic edges of the hexagons in the two
  surfaces are both equidistant curves of the same length $c/2$.

  We define the required map piecewise from each possibly truncated
  hexagon in $(Q_T,\s)$ to the corresponding one in $(Q_T,\s')$. Since
  all the Fenchel-Nielsen coordinates agree in the interior of $Q_T$,
  we only have to see that there is a bilipschitz map between the
  truncated parts of two hexagons $H$ and $H'$ with alternate
  sidelengths $l_1,l_2,l_3$ and $l_1',l_2,l_3$ coming from the pants
  curves, where $l_1,l_1' < c/2$.  Since $l_2,l_3$ are uniformly
  bounded above and below, the distance between the corresponding
  sides in both hexagons is also uniformly bounded above and below,
  see the proof of Lemma~\ref{lem:width}.  The distance between the
  side of length $l_1$ and the equidistant curve of length $c/2$ is
  equal to $\log [c/2l_1]$, up to a bounded additive error.  Hence by
  Lemma~\ref{lem:width}, the distances between the sides of lengths
  $l_2, l_3$ and the equidistant curve of length $c/2$ are bounded
  above, while they are bounded below by choice of $c$.
  
  It is now easy to define a bilipschitz homeomorphism between the
  truncations of $H$ and $H'$.  For example, fix a point $O$ whose
  distance from all sides of $H$ is uniformly bounded above and below
  and divide $H$ into six triangles by joining $O$ to the vertices of
  $H$.  Note that if we are given two hyperbolic triangles whose side
  lengths are uniformly bounded above and below, we can the map the
  three sides linearly to each other and then extend to a uniformly
  bilipschitz map on the interiors.  We can do the same even when one
  side is an equidistant curve rather than a geodesic. Now define the
  required map from $H$ to $H'$ triangle by triangle making it agree
  on the edges joining $O$ to the vertices of $H$.  It is clear that
  the resulting bound on $K$ depends only on the initial upper bound
  on the lengths of the pants curves.  Note also that $K \to 1$ as
  $\ep_0 \to 0$.
\end{proof}

We may assume that in the definition of the truncated surfaces $Q_T$,
the constants $\ep_0$ and $c$ are chosen small enough that any
non-peripheral simple geodesic loop contained in $(Q,\s)$ is
completely contained in $Q_T$. In particular, if ${\cal P}_\s$ is a
short pants system for $\s$ and if $\beta \in {\cal P}_\s \setminus
\dd Q$ is contained in $Q_T$, then so is its dual $\delta_\beta$.  Let
$M_\s$ be the short marking of $\s$ associated to ${\cal P}_\s$. We
call the subset of $M_\s$ thus defined, the {\em restriction} $M_\s
|_Q$ of $M_\s$ to $Q$ \cite{masurminskyII}. Equivalently, $M_\s |_Q$
is the set of curves in $M_\s$ that are completely contained in $Q$
and are non-peripheral in $Q$. If $Q$ is a pair of pants then $M_\s
|_Q$ is empty.

\begin{proof}[Proof of Proposition~\ref{lem:minimalonQ}]
  Where convenient, we drop the subscript $t$.  Since $Q$ is a
  component of $S \setminus {\cal S}_\ep(\mathcal G_t) $, the curves
  in $\partial Q$ are included in the set of pants curves in both
  $M_{\L}$ and $M_{\G}$.  Define a new metric $\tau=\tau_t$ on $S$
  interpolating $\mathcal G_t$ and $\L_t$ as follows.  Let $X$ be the
  metric completion of $S \setminus Q$.  First we choose a new pants
  system ${\cal P}_\tau$ for $S$.  The system ${\cal P}_\tau$ contains
  all the curves in $\dd Q$, in the interior of $Q$ it consists of the
  pants curves in $M_{\G} |_Q$, while in the interior of $X$ it
  consists of the pants curves in $M_\L |_X$.  We define $\tau_t$ by
  specifying its Fenchel-Nielsen coordinates with respect to ${\cal
    P}_\tau$. The metric $\tau_t$ will have the same Fenchel-Nielsen
  coordinates associated to the pants curves in $M_{\G}|_Q$ as
  $\mathcal G_t$ and the same Fenchel-Nielsen coordinates associated
  to the curves in $M_\L|_Q \cup \dd Q$ as $\L_t$.

Since $\L$ is on the line of minima we have:
 \begin{equation}
\label{eqn:one}
l_{\L}(\nu_t^+) + l_{\L}(\nu^-_t) 
 \leq l_{\tau}(\nu^+_t) + l_{\tau} (\nu^-_t).  
\end{equation}  
Let us estimate both sides of this inequality.  Applying
Lemma~\ref{lem:estimate} to $\nu=\nu^\pm_t$, we obtain:
$$l_\tau(\nu \cap X) = l_\tau(\nu_X) + \frac{1}{2}\sum_{\a \in
  \partial Q} l_{\tau}(\a) \, i(\a,\nu) \, Tw_{\tau}(\nu,\a) +
O(i(\nu,\partial Q)),$$
$$l_{\L}(\nu \cap X) = l_{\L}(\nu_X) + \frac{1}{2}\sum_{\a
\in \partial Q} l_{\L}(\a) \, i(\a,\nu) \, Tw_{\L}(\nu,\a) +
O(i(\nu,\partial Q)).$$
By construction, $ l_\tau(\nu_X) = l_{\L}(\nu_X)$. 
\begin{equation*} 
\begin{split}
  |\,l_\tau&(\nu \cap X) -l_{\L}(\nu \cap X)\,| \leq  \\
  & \leq \frac{1}{2}\sum_{\a \in \partial Q} \, l_{\tau}(\a) i(\nu,\a)
  \left|Tw_{\tau}(\nu,\a)-Tw_{\L}(\nu,\a) \right|+
  O\big(i(\nu,\partial Q) \big).
\end{split}
\end{equation*}
By construction, the Fenchel-Nielsen twist coordinates for $\tau$ and
$\L$ on any component $ \a$ of $\partial Q$ coincide. Therefore, by
Lemma~\ref{lem:minskytwist}, we have
$|Tw_{\tau}(\nu,\a)-Tw_{\L}(\nu,\a)| \leq 4$.  Thus
$$
|\,l_\tau(\nu \cap X)-l_{\L}(\nu \cap X)\,| \leq {2}\sum_{\a \in
  \partial Q} l_{\tau}(\a)\, i(\nu,\a) +O\big(i(\nu,\partial Q)\big).
$$
Substituting this into Equation~(\ref{eqn:one}) and noting that we
are working under the assumption that all components $ \a $ of
$\partial Q$ are extremely short in $\tau$, we obtain:
$$
l_{\L}(\nu^+_t \cap Q) + l_{\L}(\nu^-_t\cap Q) <l_{\tau}(\nu^+_t
\cap Q) + l_{\tau} (\nu^-_t \cap Q) + O(l_{q_t}(\dd Q)).
$$
Since every component of $\dd Q$ is extremely short in $\tau$, the
collar lemma implies that $i(\nu,\dd Q) \stackrel{{}_\ast}\prec
l_\tau(\nu \cap Q)$ so we may replace the last approximation by
$$
l_{\L}(\nu^+_t \cap Q) + l_{\L}(\nu^-_t\cap Q) 
\stackrel{{}_\ast}\prec l_{\tau}(\nu^+_t\cap Q) + l_{\tau} (\nu^-_t\cap Q).  
$$
Since for $\s = \L, \tau$ and $\nu=\nu_t^\pm$, we have by
Lemma~\ref{thesame}
\begin{align*}
l_\s(\nu \cap Q) 
&= l_\s(\nu \cap Q_T) + l_\s(\nu \cap (Q\setminus Q_T)) \\
& \mul l_\s(\nu_Q \cap Q_T) +l_\s(\nu \cap (Q\setminus
Q_T)),
\end{align*}
we can apply Lemma~\ref{lem:thinlength} to 
subtract the contribution of the
collars forming $Q \setminus Q_T$ from both sides to obtain:
$$
l_{\L}(\nu^+_Q \cap Q_T) + l_{\L}(\nu^-_Q \cap Q_T) 
\stackrel{{}_\ast}\prec l_{\tau}(\nu^+_Q \cap Q_T) + l_{\tau} (\nu^-_Q\cap Q_T).
$$
To complete the proof, we apply Lemmas~\ref{thesame} and 
\ref{bilipschitz}:
\begin{align*}
l_\tau(\nu^+_Q \cap Q_T) + l_\tau(\nu^-_Q \cap Q_T)
& \mul l_g(\nu^+_Q \cap Q_T) + l_{\G}(\nu^-_Q \cap Q_T) \\
&\mul l_g(\nu^+_t \cap Q_T) + l_{\G}(\nu^-_t \cap Q_T). \qedhere
\end{align*}
\end{proof}

\subsection{Correspondence between thick components}
\label{sec:correspondence}
This section contains the meat of our comparison between the
geometries of $\L_t$ and $\G_t$. We show that
(Corollary~\ref{cor:shortsame}) the sets of short curves on $\L_t$ and
$\G_t$ coincide.  Generalizing Theorem~\ref{thm:warmup}, we prove
(Theorem~\ref{thm:markingsclose} and Corollary~\ref{equiv}) that the
geometries of the thick parts of $\L_t$ and $\G_t$ are close. As in
that proof, our strategy is to use short markings to estimate lengths.
The main point is to use Proposition~\ref{lem:minimalonQ} as a
substitute for the length minimization property of $\L_t$.

\medskip

We need the following result which generalizes
Proposition~\ref{prop:capmarking} to thick components.
\begin{proposition}
\label{prop:capmarking1}
Assume that $l_\s(\a) < \ep_0$ for every component $\a$ of $ \dd Q$.
Let $\rho>0$ and suppose that $l_\s(\zeta) \geq \rho$ for every
non-peripheral simple closed curve $\zeta$ in $Q$.  Then for any
simple closed geodesic $\gamma$ on $(S,\s)$,
\begin{equation}
\label{eqn:capmarkingthick}
l_\s (\gamma \cap Q_T) \mul 
i(M_\s |_Q ,\gamma), 
\end{equation}
where the multiplicative constants depend only on $\rho$.
\end{proposition}
\begin{proof}
By Lemma~\ref{thesame} it is sufficient to prove that 
$$l_\s (\gamma_Q \cap Q_T) \stackrel{{}_\ast}{\asymp} i(M_\s |_Q
,\gamma).$$
We modify the argument in ~\cite{minskytop} Lemma 4.7.

Notice that since $Q_T$ is $\rho$-thick, it follows from
Corollary~\ref{cor:duallength} that the lengths of all the curves in
$M_\s |_Q$ are bounded above.  Cutting $Q_T$ along the curves in $M_\s
|_Q$, we obtain a collection of convex polygons $\{D_i\}$, together
with annuli $\{ A_j \}$, where one boundary component $\partial_0 A_j
$ is a component of $\partial Q_T$, while the other component
$\partial_1 A_j $ is made up of arcs in $M_\s |_Q$.

Since the total length of curves in $ M_\s |_Q$ is uniformly bounded
above, the length of $\dd D_i$ is uniformly bounded above, and
therefore, $D_i$ has uniformly bounded diameter. We claim that the
annuli $ A_j $ also have uniformly bounded diameter.  Since the length
of $\partial_0 A_j $ is bounded below by $\ep_{0}$, an area argument
shows that the distance between $\dd_1 A_j$ and $\dd_0 A_j$ is
uniformly bounded above.  Since, furthermore, the lengths of $\dd_0
A_j$ and $\dd_1 A_j$ are uniformly bounded above, it follows that $
A_j $ has uniformly bounded diameter, as claimed.  Setting
$$D = \max \{\mbox{diam} D_i, \mbox{diam} A_j\} $$ 
gives the upper bound
$$l_\s(\gamma_Q \cap Q_T) \leq i(\gamma, M_\s |_Q) \cdot D.$$

Since the lengths of all the curves in $M_\s |_Q$ are bounded above,
by the collar lemma, there is an embedded collar of definite radius
$d$ around every curve in $M_\s |_Q$. Therefore, if $\gamma$ crosses
$\beta \in M_\s |_Q$, then $l_\s(\gamma) > d\cdot i(\gamma, \beta)$.
Let $k$ be the number of pants curves in $Q$. Since there must be some
$\beta \in M_\s |_Q$ such that $i(\gamma,\beta) \geq i(\gamma,M_\s
|_Q)/(2k)$, we have
$$l_\s(\gamma_Q \cap Q_T) > \frac{d}{2k} \cdot i(\gamma,M_\s |_Q),$$
giving the desired lower bound.
\end{proof}

Applying Proposition~\ref{prop:capmarking1}, we can now deduce from
Proposition~\ref{lem:minimalonQ} that a non-peripheral curve in $Q$
cannot be too short in $\L_t$:
\begin{proposition}
\label{bddbelow}
Let $Q$ be a component of $S \setminus {\cal S}_\ep(\mathcal G_t)$
where $\ep$ is chosen as in Equation~(\ref{star}).  Then for any
non-peripheral simple closed curve $\zeta$ in $Q$, we have
$l_{\L_t}(\zeta) \stackrel{{}_\ast}{\succ} 1$.
\end{proposition}
\begin{proof}
First, we claim that 
\begin{equation}
\label{lem:comparecq}
l_{\mathcal G_t}(\nu^+_t \cap {Q_T}) + l_{\mathcal G_t}(\nu^-_t \cap {Q_T}) 
\stackrel{{}_\ast}{\asymp} \lambda_Q.
\end{equation}
To see this, let $M_{\mathcal G_t}$ be a short marking for $\mathcal
G_t$ and let $M_{\mathcal G_t} |_Q$ denote its restriction to $Q$.  By
Proposition~\ref{prop:capmarking1}, we have
\begin{align*}
l_{\mathcal G_t}(\nu^+_t \cap {Q_T}) + l_{\mathcal G_t}(\nu^-_t \cap {Q_T})
&\stackrel{{}_\ast}{\asymp} i(M_{\mathcal G_t} |_Q ,\nu^+_t)+
i(M_{\mathcal G_t} |_Q ,\nu^-_t) \\
&\stackrel{{}_\ast}{\asymp} l_{q_t}(M_{\mathcal G_t} |_Q).
\end{align*}
On the other hand,
$l_{\mathcal G_t}(M_{\mathcal G_t} |_Q) \mul 1$.  Hence by 
Theorem~\ref{thm:scalefactor}
$$
l_{q_t}(M_{\mathcal G_t} |_Q) \mul \lambda_Q$$
and the claim is proved.

Now if $\zeta$ is a non-peripheral simple closed curve in $Q$ with
$l_{\L_t}(\zeta) < \ep_0$, then consideration of the collar about
$\zeta$ gives the estimate
$$l_{\L_t}(\nu \cap {Q_T}) \stackrel{{}_\ast}{\succ}
  i(\nu, \zeta) \log \frac{1}{l_{
    \L_t}(\zeta)}$$
for any $\nu \in \ML$, so in particular,
$$l_{\L_t}(\nu^+_t \cap {Q_T}) +l_{\L_t}(\nu^-_t \cap {Q_T})
\stackrel{{}_\ast}{\succ} l_{q_t}(\zeta) \log \frac{1}{l_{
    \L_t}(\zeta)}.
$$
Proposition~\ref{lem:minimalonQ} and Equation~(\ref{lem:comparecq})
give 
$$
\lambda_Q \stackrel{{}_\ast}{\succ} l_q(\zeta) 
\log \frac 1{l_{\L_t}(\zeta)}.
$$ 
From the definition of $\lambda_Q$ we have
$\lambda_Q /l_q(\zeta) \leq 1$, so that $l_{\L_t}(\zeta)
\stackrel{{}_\ast}{\succ} 1$.
\end{proof}

\medskip

Proposition~\ref{bddbelow} and Proposition~\ref{cor:shorteronL}
together prove Theorem~\ref{shortonboth} of the Introduction, that the
sets of extremely short curves on $\L_t$ and $\mathcal G_t$ coincide.
More precisely, we can reformulate Proposition~\ref{bddbelow} as:
\begin{corollary}
\label{cor:shortsame} Let $\ep$ be as in Equation~(\ref{star}). 
Then there exists $\ep'>0$ such that $\S_{\epsilon'}(\L_t) \subset
\S_{\epsilon}(\mathcal G_t)$.
\end{corollary}
 
It is also now easy to complete the proof of our main comparison
between the thick parts of $\L_t$ and $\G_t$:
\begin{theorem}
\label{thm:markingsclose}
Let $Q$ be a component of $S\setminus {\cal S}_\ep(\mathcal G_t)$
which is not a pair of pants, and let $M_{\L_t}$ be a short marking
for $\L_t$.  Then $$l_{\mathcal G_t}( M_{\L_t} |_Q ) \mul 1.$$
\end{theorem}
\begin{proof} By Theorem~\ref{thm:scalefactor}, we have
\begin{align*}
  l_{\mathcal G_t}(M_{\L_t} |_Q)
  &\mul \frac{1}{\lambda_Q}l_{q_t}(M_{\L_t} |_Q)\\
  &\mul \frac{1}{\lambda_Q} [ i(M_{\L_t} |_Q, \np) + i(M_{\L_t}
  |_Q,\nm)] .
\end{align*}  
By Proposition~\ref{bddbelow}, there is a constant
$\rho=\rho(\ep_0)$ depending only on $\ep_0$ such that
$l_{\L_t}(\zeta) > \rho(\ep_0)$ for every non-peripheral curve
$\zeta$ in $Q$. Therefore, we can
apply Proposition~\ref{prop:capmarking1} to get
$$
i(M_{\L_t} |_Q, \nu^+_t) + i(M_{\L_t} |_Q,\nu^-_t) \mul
l_{\L_t}(\nu^+_t \cap Q_T) + l_{\L_t}(\nu^-_t \cap Q_T).$$
Since the
lower bound $l_{\mathcal G_t}(M_{\L_t} |_Q) \stackrel{{}_\ast}{\succ} 1$ is
trivial, the result now follows from Proposition~\ref{lem:minimalonQ}
and Equation~(\ref{lem:comparecq}).
\end{proof} 
Equivalently, we can formulate Theorem~\ref{thm:markingsclose} in
terms of the surface $Q_0$ obtained from $Q$ by replacing every
boundary component with a puncture.  Let $\L_t |_{Q_0}$ and $\mathcal
G_t |_{Q_0} $ be respectively, the surface $Q_0$ equipped with the
metrics obtained from $\L_t$ and $\mathcal G_t$ by pinching the curves
in $\dd Q$ but otherwise leaving the metric unchanged. In other words,
in the notation of the product region theorem in
Section~\ref{sec:productregion}, the collection of curves in $\dd Q$
is $\mathcal A$ and the metrics on $Q_0$ are defined by $\Pi_0(\L_t)$
and $\Pi_0(\G_t)$, respectively, restricted to the component $Q_0$ of
$S_{\mathcal A}$.  Then we have:
\begin{corollary}
\label{equiv}
Let $Q$ be a component of $S\setminus {\cal S}_\ep(\mathcal G_t)$. Then
  $$d_{{\cal T} (Q_0)} (\L_t |_{Q_0},\mathcal G_t |_{Q_0}) = O(1).$$
\end{corollary}
\begin{proof}
  The boundary components of both $(Q,\L_t)$ and $(Q,\mathcal G_t)$
  are extremely short. In this case, it was shown in
  \cite{minskyproduct} (see the proof of Lemma~\ref{bilipschitz}) that
  $(Q,\L_t)$ and $(Q,\mathcal G_t)$ can be embedded into $(Q_0,\L_t
  |_{Q_0})$ and $(Q_0,\mathcal G_t |_{Q_0})$ respectively, by a
  $K$-quasi-conformal map, where $K$ depends only on $\ep_0$.  Since
  simple curves do not penetrate the thin part of $Q_0$, the
  restriction $M_{\L_t}|_Q$ is a short marking for $(Q_0,\L_t
  |_{Q_0})$. By Theorem~\ref{thm:markingsclose}, we have $l_{\mathcal
    G_t}(M_{\L_t}|_Q) \mul 1$. The result now follows as in the proof
  of Theorem~\ref{thm:warmup}.
\end{proof}
\subsection{Comparison of lengths of short curves}
\label{sec:compshtlength}
Theorem~\ref{thm:markingsclose} allows us to complete the proof of
Equation~(\ref{HisK}): \begin{proposition}
\label{lem:HandK} Let $\a$ be an extremely short curve on $\L_t$. Then 
 $$H_t(\a)  \prec K_t(\a). $$
\end{proposition} 
\begin{proof}
  With $\cal B$ as in Equation~(\ref{eqn:defofH}), let $\beta \in \cal
  B$ be the curve that has the largest $q_t$-length, so that
$$H_t(\a) \mul \frac{l_{q_t}(\beta)}{l_{q_t}(\a)}.$$

Since ${\cal S}_\ep(\mathcal G_t) \subset {\cal S}_{\ep_0}(\L_t)$, it
follows that the curves in $\cal P_{\L_t}$ are disjoint from $\cal
S_\ep(\mathcal G_t)$.  Thus, $\a$ and $\beta$ are contained in the
closure of a common component $Q$ of $S \setminus \cal S_\ep(\mathcal
G_t)$.

Suppose first that $\beta$ is not peripheral in $Q$.  Then $\beta \in
M_{\L_t}|_Q$, so that by Theorem~\ref{thm:scalefactor}
and~\ref{thm:markingsclose},
$$l_{q_t}(\beta) \mul \lambda_Q l_{\mathcal G_t}(\beta) \mul \lambda_Q.$$
If in
addition, $\a$ is not peripheral, then $l_{q_t}(\a) \mul \lambda_Q$ so
that $H_t(\a) \mul 1$ and the desired inequality holds trivially. If
$\a$ is peripheral, then
$$H_t(\a) \mul  \frac{l_{q_t}(\beta)}{l_{q_t}(\a) }\mul
\frac{\lambda_Q}{l_{q_t}(\a)} \leq K_t(\a). 
$$

Now suppose that $\beta$ is peripheral in $Q$. If $Q$ is a pair of
pants, then the desired inequality follows from the definition of
$H_t(\a)$ and $K_t(\a)$. If $Q$ is not a pair of pants, then since the
component of $Q \setminus M_{\mathcal G_t}|_Q$ containing $\beta$ is
an annulus whose one boundary component is $\beta$ and the other a
finite (at most $4$) union of arcs coming from curves $\cup \gamma_i$
in $M_{\mathcal G_t} |_Q$, again by Theorem~\ref{thm:scalefactor} we
obtain
$$l_{q_t}(\beta) \leq \sum l_{q_t}(\gamma_i)
\mul \sum \lambda_Q l_{\mathcal G_t}(\gamma_i) \mul \lambda_Q,$$   
from which the result follows as before.
\end{proof}

Theorem~\ref{shortonL} now follows immediately from
Theorem~\ref{thm:shortonL}, Proposition~\ref{cor:shorteronL}, and
Proposition~\ref{lem:HandK}, completing our comparison between short
curves on $\G_t$ and $\L_t$.
 \begin{theorem} 
\label{cor:final}
Let $\a$ be any curve on $S$ which is neither vertical nor horizontal.
If $\a$ is extremely short in $\L_t$, then
$$\frac{1}{l_{\L_t}(\a)} \asymp \max \left\{ D_t(\a), \sqrt{K_t(\a)}
\right\}.$$
\end{theorem}

In case $\a$ is vertical or horizontal, we have
\begin{theorem} 
\label{cor:final1}
If $\a$ is vertical, then
$$
\frac{1}{l_{\L_t}(\a)} 
\asymp \max \left\{e^{-2t} \mod F_0(\a), \sqrt{K_t(\a)} \right\}.
$$
If $\a$ is horizontal, then the estimate is the same except
that the first term is replaced by $e^{2t} \mod F_0(\a).$
\end{theorem}
\begin{proof}
  There are multiplicative constants depending only on the fixed
  laminations $\nu^\pm$ such that
  $$\frac{m^\mp(\a)}{i(\nu^\pm,\a)} \mul \mod F_0(\a)$$
  (see
  Theorems~\ref{thm:shortonL1} and~\ref{prop:shortingeo1}) holds
  independently of $\a$ in a tautological way, due to the fact that
  the total number of vertical (or horizontal) curves is finite; it is
  bounded above by $-\chi(S)$.  The proofs of
  Proposition~\ref{cor:shorteronL} and~\ref{lem:HandK} go through in
  this case, so that $H_t(\a) \asymp K_t(\a)$.  Hence the estimate
  follows from Theorem~\ref{thm:shortonL1}.
\end{proof}

\subsection{\Teich distance}
With the preceding collection of results at our disposal,
Theorem~\ref{distance} of the Introduction becomes an easy application
of Minsky's product region theorem~\ref{thm:minskyproduct}.
\label{sec:finalcomparison}
\begin{theorem}
 \label{thm:finalcomparison} 
 The \Teich distance between $\G_t$ and $\L_t$ is given by
 $$
 d_{\ts} (\G_t, \L_t) = \max_{\a \in {\cal S}_\ep(\mathcal G_t)}
 \frac{1}{2} \Bigg| 
\log \frac{l_{\mathcal G_t}(\a)}{l_{\L_t}(\a)} \Bigg| \pm O(1).$$
\end{theorem}
\begin{proof} 
  As noted before, $l_{\L_t}(\a) \leq \ep_0$ for every $\a \in {\cal
    S}_\ep(\mathcal G_t)$. To simplify notation, let $\E_t = {\cal
    S}_\ep(\mathcal G_t)$.  By Theorem~\ref{thm:minskyproduct}, we
  have
\begin{equation*}
\begin{split}
d_{\ts} &(\G_t, \L_t) = \\
&= \max_{\a \in \E_t} \{
  d_{\T(S_{\E_t})}(\Pi_0(\G_t),\Pi_0(\L_t)),
  d_{\HH_{\a}}(\Pi_{\a}(\G_t),\Pi_{\a}(\L_t)) \} \pm O(1),
\end{split}
\end{equation*}
where $S_{\E_t}$ is the surface obtained from $S$ by removing $\E_t$
and replacing the resulting boundary components by punctures and
$\Pi_0, \Pi_\a$ are defined as in Section~\ref{sec:productregion}.
From Corollary~\ref{equiv} we deduce that $$
d_{\T(S_{{\cal
      E}})}(\Pi_0(\G_t),\Pi_0(\L_t)) = O(1).$$
If $t-t_\a >0$, then by
Theorem~\ref{prop:quadtwistest} or~\ref{prop:quadtwistest1} we have
$$Tw_{\mathcal G_t}(\nu^+,\a)l_{\mathcal G_t}(\a) = O(1)$$
and by
Theorem~\ref{thm:twistonL} or~\ref{thm:twistonL1} we have
$$Tw_{\L_t}(\nu^+,\a)l_{\L_t}(\a) = O(1).$$  
Applying Corollary~\ref{cor:hypcomparison} to $\mathcal G_t$ and $\L_t$ and 
the lamination $\nu^+$, we find
$$\exp{2d_{\HH_{\a}}(\Pi_{\a}(\G_t),\Pi_{\a}(\L_t))} \mul
l_{\mathcal G_t}(\a)/l_{\L_t}(\a) .$$
If $t- t_{\a} <0$, the same result holds
by applying a similar argument with $\nu^-$.
\end{proof}

\subsection{Examples}
The combinatorial nature of our length estimates allows us to use
Theorem~\ref{thm:finalcomparison} to construct examples in which $\cal
L$ and $\cal G$ have a variety of different relative behaviors. As a
special case, if $S$ is a once-punctured torus or four-times-punctured
sphere, every thick component must be a pair of pants. In this case,
$K_t(\a)$ is bounded and therefore, a curve gets short only if
$d_\a(\np,\nm)$ is large:
\begin{corollary} If $S$ is a once-punctured torus or a four-times
  punctured sphere, then for any measured laminations $\np,\nm$, the
  associated \Teich geodesic and line of minima satisfies
$$d_{\ts}(\G_t,\L_t) = O(1).$$
\end{corollary}

\medskip On surfaces of higher genus, it is possible to have $\np,\nm$
and $\a$ such that $d_\a(\np,\nm)$ is bounded while $K_t(\a)$ is
arbitrarily large. We can construct a simple example as follows. Take
two Euclidean squares each of area $1/2$ and cut open a slit of length
$\varepsilon$ at each of their centers. Although it is not necessary,
for concreteness we can assume that in both squares, the slit is
parallel to a pair of sides.
\begin{figure}[htb]
\centerline{\epsfbox{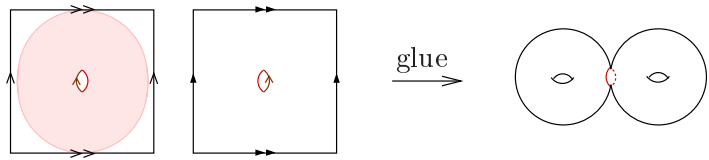}}
\caption{Glue two slit tori along slit.} 
\label{fig:example}
\end{figure}
Foliate each square by the two mutually orthogonal foliations that
both make angle $\pi/4$ with the slit and for each, take the
transverse measure induced by the Euclidean metric. Identify pairs of
sides in each square to obtain two one-holed tori $T_1,T_2$ and glue
$T_1,T_2$ along their boundaries, as shown in the figure, to obtain a
genus two surface $S$ with waist curve $\a$.  The two foliations match
along $\a$ and specifying one to be the vertical foliation defines a
quadratic differential $q=q_0$ on $S$. Let $\G$ and $\L$ be
respectively, the \Teich geodesic and the line of minima defined by
the vertical and horizontal foliations $\nm$ and $\np$ of $q_0$. In
this example, the foliations are rational, but it is easy to see that
varying the initial angle of the slit gives more general foliations.

Let $q_t$ be the associated family of quadratic differentials.  Note
that at time $t=0$ the curve $\a$ is balanced.  The $q_t$--geodesic
representative of $\alpha$ is unique and the flat annulus
corresponding to $\alpha$ is degenerate. Thus, by
Proposition~\ref{prop:flat}, we have $d_\alpha(\np, \nm) =O(1)$.  But
since $\lambda_{T_1}=\lambda_{T_2} =1/\sqrt{2}$ at $t=0$, we have
$K_0(\alpha) \asymp 1/\varepsilon$ (the shaded region indicates a
maximal expanding annulus).  Assuming $\varepsilon$ is very small,
Theorems~\ref{prop:shortingeo} and~\ref{thm:shortonL} give
$$
\frac{1}{l_{\G_0}(\a)} \asymp \log \frac{1}{\varepsilon} 
\quad\text{and}\quad
\frac{1}{l_{\L_0}(\a)} \asymp \frac{1}{\sqrt{\varepsilon}}.
$$
Thus by Theorem~\ref{thm:finalcomparison},
$$d_{\ts}(\G_0,\L_0) 
\succ \frac{1}{2} \log \frac{l_{\G_0}(\alpha)}{l_{\L_{0}}(\alpha)}
\asymp \log \frac{1}{\sqrt{\varepsilon}
\log[1/\varepsilon]}.$$

In fact, because for any two hyperbolic metrics $\s,\tau$ we
have \cite{wolpert2}
$$d_{\ts}(\s,\tau) \geq \frac{1}{2}\log \sup_{\zeta \in \cal S}
\frac{l_\s(\zeta)}{l_\tau(\zeta)}$$ and because 
the length of $\alpha$ along $\G$ is (coarsely)
shortest at the balance time $t_\a = 0$ \cite{rafi1}, 
the following stronger inequality holds:
$$
\inf_{t \in \RR} d_{\ts}(\G_t, \L_0) \geq 
\inf_{t \in \RR} 
\frac{1}{2}\log \frac{l_{\mathcal G_t}(\alpha)}{l_{\L_{0}}(\alpha)} \succ
\frac{1}{2}\log \frac{l_{\G_0}(\alpha)}{l_{\L_{0}}(\alpha)}.
$$
Taking $\varepsilon$ small enough we can ensure that $\L_0 $ is as
far as we like from any point on $\G$. This example can be easily
extended to any surface of large complexity, by which we mean a
surface whose genus $g$ and number of punctures $p$ satisfies $3g-4 +
p \geq 1$. In summary,
\begin{corollary} 
  If $S$ is a surface of large complexity, then given any $n>0$, there
  are measured laminations $\nu^+(n), \nu^-(n)$ on $S$ which depend on
  $n$, such that for the associated \Teich geodesic $\G(n)$ and line
  of minima $\L(n)$,
  $$\inf_{t \in \RR} d_{\ts}(\G_t(n),\L_0(n))  > n.$$
\end{corollary}
It is also possible to construct examples for any such surface where
the two measured laminations are fixed and the associated \Teich
geodesic and line of minima satisfy $d_{\ts}(\G_{t_n},\L_{t_n}) > n$
for a sequence of times $t_n \to \infty$ as $n \to \infty$. This
however, is beyond the scope of this paper.
\section{Appendix}
We give proofs of the length estimates that were
deferred in previous sections.\\

\begin{figure}[htb]
\centerline{\epsfbox{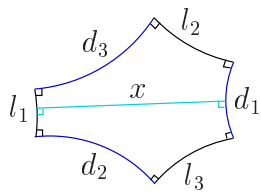}}
\caption{Half a pair of pants.} 
\label{fig:hex1}
\end{figure}

\noindent{\em Proof of Lemma~\ref{lem:width}.}
Let $H$ be one of the two right-angled hexagons obtained by cutting
$P$ along its three seams.
Let $l_i = l(\a_i)/2$ and $d_i$ be the length
of the seams, labeled as shown in
Figure~\ref{fig:hex1}.
By the cosine formula for right-angled hexagons, we have
$$\cosh d_3 = \frac{\cosh l_3 + \cosh l_1 \cosh l_2}{\sinh l_1 \sinh
  l_2}.$$
By hypothesis, $l_i < L/2$, so $\sinh l_i
\mul l_i$ and $\cosh l_i \mul 1$,
where the multiplicative constants involved depend only on $L$.
Therefore, 
$$
\cosh d_3  \mul  \frac{1}{l_1l_2}.
$$
It follows that $d_3 = \log[1/l_1] + \log[1/l_2] \pm O(1)$, where the
bound on the additive error depends only on $L$.

Now we estimate the length of the perpendicular from $\a_i$ to itself.
Let $x$ be the length of the perpendicular as in Figure~\ref{fig:hex1}.
By the formula for right-angled pentagons, we have
$$\cosh x = \sinh l_2 \sinh d_3.$$ Since $\sinh l_2 \mul l_2$ and by
the above, $\sinh d_3 \mul 1/[l_1 l_2]$, it follows that $\cosh x \mul
1/l_1$.  Hence, $x= \log [1/l_1] \pm O(1).$ Since $P$ is made of two
isometric copies of $H$, we obtain the desired estimate.
\hfill $\Box$\\

 To prove Lemma~\ref{lem:hexagon*},
in addition to the standard hexagon and pentagon formulae (see for
example~\cite{beardon} or~\cite{cmswolpert}), we need the following
expression for derivatives of side-lengths derived
in~\cite{cmswolpert} Proposition 2.3:
\begin{lemma} Let $H$ be a planar right-angled hexagon with sides 
labeled $i=1,\ldots,6$ in cyclic order about $\partial H$. Let 
 $l_i $ denote the  length of side $i$ and for $n$ mod $6$, let $p_{n,
n+3}$ be the perpendicular distance from side $n$ to side $n+3 $.
Letting ${'}$ denote  derivative with respect to some variable $x$, we have
 \begin{equation}
 \label{eqn:hexderiv}
 ({\cosh p_{n,n+3}}) l_n^{'} = l_{n+3}^{'} - (\cosh{l_{n-2}}) l_{n-1}^{'}
-(\cosh{l_{n+2}}) l_{n+1}^{'}.
 \end{equation}
\end{lemma}

It is convenient to subdivide Lemma~\ref{lem:hexagon*} into two parts,
Lemma~\ref{lem:hexagon} and Lemma~\ref{lem:proofperp}, depending on
whether or not the common perpendicular in question is adjacent to
$\a$. We begin with a somewhat more detailed discussion of the
possible configurations.

\begin{figure}[htb]
\centerline{\epsfbox{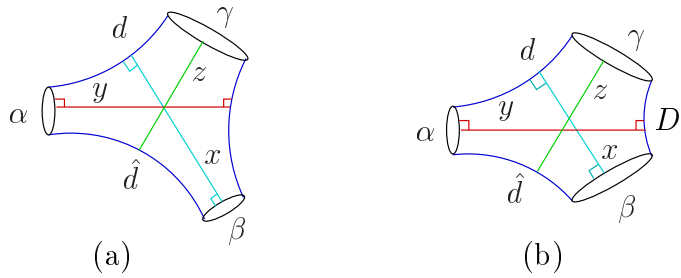}}
\caption{(a) $\pi(\a)=\pi(\beta)$;\ \ (b) $\pi(\a) \neq \pi(\beta)$.}
\label{fig:perp1}
\end{figure}

Let $P$ be a pair of pants in $S\setminus{\cal P_{\L_t}}$ that has
$\a$ as a boundary component.  For clarity, we distinguish between the
three boundary curves $\a,\beta,\gamma$ of $P$ and their projections
$\pi(\a),\pi(\beta),\pi(\gamma)$ to $S$. We may always assume that
$\pi(\gamma) \neq \pi(\a)$. There are then two possible configurations
depending on whether or not $\pi(\beta) = \pi(\a)$.
Figure~\ref{fig:perp1}(a) represents the case in which $\pi(\beta) =
\pi(\a)$ and Figure~\ref{fig:perp1}(b), the case in which $\pi(\beta)
\neq \pi(\a)$.  In (b), we do not rule out the possibility that
$\pi(\beta)= \pi(\gamma)$.  This leads to a dichotomy in the formulae
used in the proofs, but not in the final estimates.  Let $d$ be the
length of the perpendicular between $\a$ and $\gamma$ and let $h_\a$
be the length of the shortest perpendicular from $\a$ to itself.

\begin{lemma}
\label{lem:hexagon} Suppose that $\a$ is extremely short in $\L_t$. 
Then the derivatives of the perpendiculars adjacent to $\a$ are as
follows:
$$
\text{(i)} \ d'=\frac {\dd d}{\dd  l(\a)} \mul -\frac{1}{l(\a)}, \qquad
\text{(ii)} \ h_\a'=\frac{\dd h_\a}{\dd l(\a)}  \mul -\frac{1}{l(\a)}. 
$$
\end{lemma}
\begin{proof}
  Let $x,y,z$ be the lengths of
the perpendiculars as shown in Figure~\ref{fig:perp1}. 

\noindent (i) In case (a), the formula~(\ref{eqn:hexderiv}) and the pentagon
formula (\cite{cmswolpert} lemma 2.1) give, respectively,
$$d'\cosh x = - \frac{\cosh \hat d}{2}, \ \ \cosh x = \sinh \hat d
\sinh \frac{l(\a)}{2}.$$
If $l(\a)$ is small, then $\sinh l(\a) \mul
l(\a)$ and by Lemma~\ref{lem:width}, $\hat d \succ 1/\log l(\a)$ so
that $\coth \hat d = O(1)$. Hence
$$d' = -\frac{\coth \hat d}{\sinh [l(\a)/2]} \mul -\frac{1}{l(\a)}.$$

In case (b), using the same formulae as above, we get
$$d'\cosh x= \frac{1 - \cosh \hat d}{2}, \ \ \cosh x = \sinh \hat d
\sinh \frac{l(\a)}{2}. $$ 
Now by Lemma~\ref{lem:width}, $\hat d = 2 \log [1/l(\a)] \pm O(1)$ and
therefore
$\sinh \hat d \mul 1/l(\a)^2 \mul \cosh \hat d.$ Hence,
$$d' = \frac{1 - \cosh \hat d}{2 \sinh \hat d \sinh [l(\a)/2]}
\mul l(\a) \left[1-\cosh \hat d \right] \mul -\frac{1}{l(\a)}.$$
(ii) In case (a), $h_\a = 2y$ and $\cosh y = \sinh d
\sinh[l(\gamma)/2]$. Differentiating both sides, we get 
$$y'= \frac{\sinh[l(\gamma)/2] \cosh d}{\sinh y} \ d'.$$ 
By Lemma~\ref{lem:width},
$\cosh d \mul 1/[l(\a)l(\gamma)]$ and
$\sinh y \mul 1/l(\a)$, so we obtain $y' \mul -1/l(\a)$, as
desired.

In case (b), $h_\a= \hat d$ and by \cite{cmswolpert} equation(6),
$$\hat d' \cosh z = - \cosh d.$$ Substituting 
$\cosh z \mul 1/l(\gamma)$ and $\cosh d \mul 1/[l(\a)l(\gamma)]$
gives the desired estimate.
\end{proof}

\begin{figure}[htb]
\centerline{\epsfbox{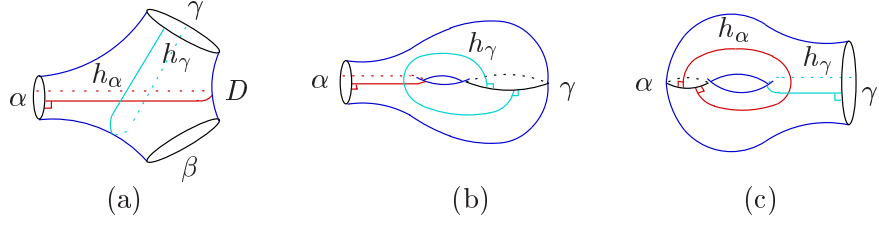}}
\caption{(a) $\pi(\a),\pi(\beta),\pi(\gamma)$ all distinct; 
(b) $\pi(\beta)=\pi(\gamma)$; (c) $\pi(\a)= \pi(\beta)$}
\label{fig:perp2}               
\end{figure}

Now consider perpendiculars in $P$ that are disjoint from $\a$.
Let $h_\gamma$ be the length of the perpendicular from $\gamma$ to
itself, as shown in Figure~\ref{fig:perp2}.  Further, if $\pi(\beta)
\neq \pi(\a)$, let $D$ be the length of the perpendicular between
$\beta,\gamma$. We have $h_\gamma = D$ when $\pi(\beta) = \pi(\gamma)$
(see Figure~\ref{fig:perp2}(b)).
\begin{lemma}
\label{lem:proofperp}
Suppose that $\a$ is extremely short in $\L_t$. Then the derivatives
of the perpendiculars not adjacent to $\a$ are as follows:
$$
\text{(i)} \ h'_\gamma = \frac{\dd h_\gamma}{\dd l(\a)} 
\mul l(\a), \qquad
\text{(ii)}\ D'=\frac{ \dd D}{\dd l(\a)}  
      \mul l(\a).
$$
\end{lemma}
\begin{proof}
  (i) Assume that $\pi(\beta) \neq \pi(\gamma)$ so that we are in the
  configuration of Figure~\ref{fig:perp2} (a) or (c).  Consider the
  `front' hexagon in $P$ and denote the lengths of the sides as shown
  in Figure~\ref{fig:hex}, so that $l_1=l(\a)/2$, $l_2 = l(\gamma)/2$,
  and $z=h_\gamma/2$.
\begin{figure}[htb]
\centerline{\epsfbox{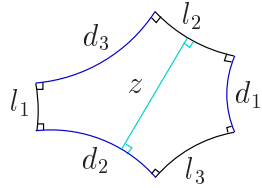}}
\caption{Front right-angled hexagon.}
\label{fig:hex}
\end{figure}
By the pentagon formula, $\cosh z = \sinh d_3 \sinh l_1.$ Taking the
derivative with respect to $l_1$, we get
$$z'\sinh z = d_3'\cosh d_3 \sinh l_1 + \sinh d_3 \cosh l_1.$$
In the
proof of Lemma~\ref{lem:hexagon}, we had $d_3'=-\coth d_2/\sinh l_1$,
and by the cosine formula for right-angled hexagons, we have
$$\cosh l_1 = \frac{\cosh d_1 + \cosh d_2 \cosh d_3}{\sinh d_2 \sinh
  d_3}.$$
Substituting these, we get
$$z'=\frac{1}{\sinh z} \cdot \frac{\cosh d_1}{\sinh d_2}.$$
Now by the
sine formula for right-angled hexagons,
$$\frac{1}{\sinh d_2} = \frac{\sinh l_1}{\sinh d_1 \sinh l_2}.$$
With this, we have
$$
z'=\sinh l_1 \cdot \frac{\cosh d_1}{\sinh d_1} \cdot
\frac{1}{\sinh l_2 \sinh z} 
\mul
\sinh l_1 \mul l_1,$$
since $\coth d_1 \mul 1$ and $\sinh z \mul 1/l_2$
when $l_1,l_2$ are respectively, bounded.

In the case $\pi(\beta)=\pi(\gamma)$, we have $h_\gamma' = D'$,
which is computed below.\\ 
(ii) Let $y$ be the length of the perpendicular between $\a$ and the
common perpendicular of $\beta,\gamma$, as in
Figure~\ref{fig:perp1}. Then by Equation~(\ref{eqn:hexderiv})
and the pentagon formula, we have
$$D' \cosh y= 1,\ \ \cosh y=\sinh d \sinh \frac{l(\gamma)}{2}.$$
By Lemma~\ref{lem:width}, 
$$d = \log [1/l(\a)] + \log [1/l(\gamma)] \pm O(1)$$ 
and therefore $\sinh d \mul 1/[l(\a) l(\gamma)].$ Since the
pants system is short, $\sinh [l(\gamma)/2] \mul l(\gamma)$. Thus $D'
\mul l(\a)$, as claimed.
\end{proof}

Lemmas~\ref{lem:hexagon} and~\ref{lem:proofperp} together prove
Lemma~\ref{lem:hexagon*}.

\medskip

\begin{proof}[Proof of Lemma~\ref{lem:estimate}]
If $\xi$ has a component $\varphi$ whose support is contained in $Q$,
then $\varphi \cap Q = \varphi = \varphi_Q$, which has no effect on the
inequality. Thus, we assume that no component of $\xi$ has support
entirely contained in $Q$.  Then, for simplicity, let us further assume that
$\xi$ is a simple closed curve.  Since both sides of the inequality
scale linearly with weights, it is sufficient to prove the lemma under
this assumption.  The basic idea is to approximate an arc $\eta$ of
$\xi \cap Q$ with the union of $\eta_Q$ and segments $p\hat p$, $q\hat
q$ which run along $\dd P$ between the endpoints $\hat{p},\hat{q}$ of
$\eta_Q$ and the corresponding endpoints $p,q$ of $\eta$. Let
$\a_p,\a_{q}$ denote the components of $\dd P$ that contain $p,q$,
respectively.  It is possible that $\a_p = \a_q$.

We will show that there are uniform constants $C,C'$ such that
\begin{eqnarray}
&&|\, l_\s(p\hat p) - l_\s(\a_p)\cdot 
Tw_\s(\xi,\a_p)/2 \,| < C \label{eqn:segment}\\
&&
|\, l_\s(q\hat q) - l_\s(\a_q)\cdot 
Tw_\s(\xi,\a_q)/2 \,|<C \nonumber \\
&&
|\,  l_\s(\eta) -[l_\s(p\hat p) + l_\s(\eta_Q) 
+ l_\s (q\hat q) ]\,|<C'.
\label{eqn:discrepancy}
\end{eqnarray}
It is convenient to consider the picture in the universal cover
$\HH^2$, as shown in Figure~\ref{fig:case1shrunk}.
\begin{figure}[htb]
\centerline{\epsfbox{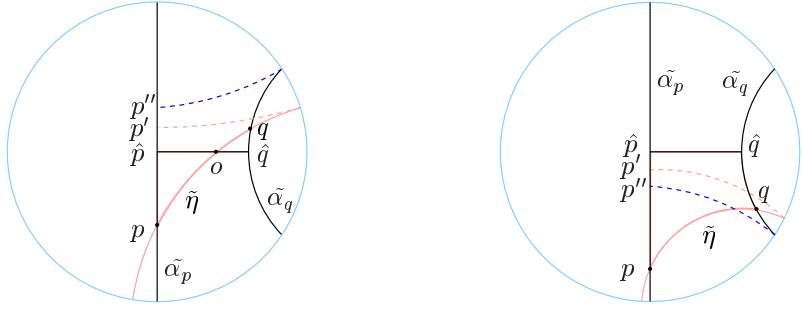}}
\caption{Approximating the length of $\tilde \eta$.} 
\label{fig:case1shrunk}
\end{figure}
Let $\tilde \eta$ be a lift of $\eta$ and let
$\tilde{\a_p},\tilde{\a_q}$ be the lifts of
$\a_p,\a_q$ that contain the endpoints of $\tilde \eta$.
Since $\eta_Q$ is homotopic to $\eta$ relative to $\dd Q$ and is
perpendicular to $\dd Q$, its lift $\tilde \eta_Q$ is the unique
perpendicular between $\tilde{\a_p},\tilde{\a_q}$, drawn
as the segment $\hat p \hat q$ in the figure. There are
two possible cases, depending on whether or not $\tilde \eta$
intersects $\tilde \eta_Q$.

Let $p',p''$ be the feet of the perpendiculars as shown. That is, $p'$
is the foot of the projection from the geodesic ray extending $\tilde
\eta$ to $\tilde \a_p$ and $p''$ is the foot of the projection
from $\tilde \a_q$ to $\tilde \a_p$. In either case,
 by definition of twist,
$$Tw_\s(\xi,\a_p,p)=2 l_\s(pp')/ l_\s(\a_p)$$ 
and furthermore,
\begin{equation}
\label{eqn:presegment}
|\,l_\s(p\hat p) - l_\s(pp') \,| < l_\s(\hat p p'').
\end{equation}
On the other hand, by trigonometry in $\HH^2$, we have
$$\cosh l_\s(\hat p p'')= 1/\tanh l_\s(\hat p \hat q) =1/\tanh l_\s
(\eta_Q).$$
Since $\eta_Q$ goes from $\a_p$ to $\a_q$ and since by
hypothesis $l_\s(\a_p),l_\s(\a_q) < \ell$, it follows from the collar
lemma that $l_\s(\eta_Q)>c(\ell)$ for some constant $c(\ell)$
depending only on $\ell$. This implies that there is a constant
$C=C(\ell)$ depending only on $\ell$ such that $l_\s(\hat p p'') <
C(\ell).$ Therefore, Equation~(\ref{eqn:presegment}) gives
Equation~(\ref{eqn:segment}), as desired. Of course, the same argument
applies to $\a_q$ so
\begin{equation}
\label{eqn:presegmentq}
|\, l_\s(q\hat q) - l_\s(\a_q)\cdot 
Tw_\s(\xi,\a_q,q)/2 \,| <C.
\end{equation}

To show Equation~(\ref{eqn:discrepancy}), we use the well known fact
that for any $\theta_0>0$, there exists a constant $k(\theta_0)$ such
that for any hyperbolic triangle with sidelengths $a,b,c$ and angle
$\theta$ opposite to $c$ with $\theta \geq \theta_0$, we have $a+b -c
<k(\theta_0)$. In the case where $\tilde \eta$ intersects $\tilde
\eta_Q$, as in the figure on the left, we apply this to the triangles
$\triangle op\hat p$ and $\triangle oq\hat q$ and get
$$ l_\s(p\hat p) + l_\s(\hat p \hat q) + l_\s(q\hat q)-l_\s(pq)< k(\pi/2).$$

In the case on the right, we apply this to triangles $\triangle p\hat
q \hat p$ and $\triangle pq\hat q$. To see that $\angle q\hat q p$ is
bounded below by some $\theta_0$, observe that 
$$\angle q\hat q p = \pi/2 - \angle p\hat q \hat p$$ 
and that 
$$\sin(\angle p\hat q \hat p)< \frac1{\cosh l_\s(\hat p \hat q)}.$$ 
Since $l_\s(\hat p \hat q) =l_\s(\tilde
\eta_Q) > c(\ell)$, it follows that $\angle p\hat q \hat p$ is bounded
away from $\pi/2$ and so $\angle q\hat q p$ is bounded below by some
constant $\theta_0=\theta_0(\ell)$, as desired. Thus in this case,
$$ l_\s(p\hat p) + l_\s(\hat p \hat q) + l_\s(q\hat q)
-l_\s(pq) <k(\pi/2) + k(\theta_0(\ell)),$$
completing the proof of Equation~(\ref{eqn:discrepancy}).  

Combining  Equations~(\ref{eqn:segment}),(\ref{eqn:discrepancy}), and
(\ref{eqn:presegmentq}) we obtain
$$\Bigg|\, l_\s(\eta)-\bigg[ l_\s(\eta_Q) +  l_\s(\a_q) 
\frac{Tw_\s(\xi,\a_q)}{2}  + l_\s(\a_p)
\frac{Tw_\s(\xi,\a_p)}{2} \bigg]\Bigg| < K(\ell).$$
Summing over all arcs $\eta$ in $\xi \cap Q$, we obtain
\begin{equation*}
\begin{split}
\Bigg| \, l_\s(\xi \cap Q)-\bigg[ l_\s(\xi_Q) + \sum_j l_\s(\a_j)
\frac{Tw_\s(\xi,\a_j)}{2}  &  i(\xi,\a_j) \bigg] \Bigg| < \\
&< K(\ell) \, i(\xi,\dd Q).  \qedhere
\end{split} 
\end{equation*}
\end{proof}

\medskip

\end{document}